                      \def\version{13 December 2005}                            %

\documentclass[reqno,11pt]{amsart} 
\usepackage{amsmath} 
\usepackage{amssymb} 
 
\newcommand{\sca}{{\alpha}}   
\newcommand{\bdot}{\,{\boldsymbol\cdot}\,}
\def\bar{\begin{array}} 
\def\ear{\end{array}} 
\def\be{\begin{equation}} 
\def\ee{\end{equation}} 
\def\bear{\begin{eqnarray*}} 
\def\eea{\end{eqnarray*}} 
\def\beal{\begin{eqnarray}} 
\def\eeal{\end{eqnarray}} 
\def\bit{\begin{itemize}} 
\def\eit{\end{itemize}} 
 
\def\d{\delta}

\def\m{\mu}

\newfam\Bbbfam 
\font\tenBbb=msbm10 
\font\sevenBbb=msbm7 
\font\fiveBbb=msbm5 
\textfont\Bbbfam=\tenBbb 
\scriptfont\Bbbfam=\sevenBbb 
\scriptscriptfont\Bbbfam=\fiveBbb

\newcommand{\R}     {\mathbb{R}} 
 
\newcommand{\N}     {\mathbb{N}} 
\renewcommand{\P}   {\mathbb{P}} 
 
\newcommand{\E}     {\mathbb{E}} 
 \newcommand{\floor}[1]{\left\lfloor #1 \right\rfloor}

\def\1{{\mathchoice {1\mskip-4mu\mathrm l}      
{1\mskip-4mu\mathrm l} 
{1\mskip-4.5mu\mathrm l} {1\mskip-5mu\mathrm l}}} 
\newcommand{\ssup}[1] {{\scriptscriptstyle{({#1}})}} 
\def\comment#1{} 
\newtheoremstyle{thm}{2ex}{2ex}{\itshape\rmfamily}{} 
{\bfseries\rmfamily}{}{1.7ex}{} 
 
\newtheoremstyle{rem}{1.3ex}{1.3ex}{\rmfamily}{} 
{\itshape\rmfamily}{}{1.5ex}{} 
 
\newenvironment{proofsect}[1] 
{\vskip0.1cm\noindent{\bf #1.}\hskip0.5cm}

\newtheorem{theorem}{Theorem}[section] 
\newtheorem{lemma}[theorem]{Lemma} 
\newtheorem{prop}[theorem] {Proposition} 
\newtheorem{cor}[theorem]  {Corollary}

\newtheorem{step}{STEP} 
 
\newcommand{\en}       {\end{equation}} 
\newcommand{\eq}       {\begin{equation}} 
 
\newcommand{\eqry}   {\begin{eqnarray}} 
\newcommand{\enqry}   {\end{eqnarray}} 
\newcommand{\eqarray}   {\begin{eqnarray}} 
\newcommand{\enarray}   {\end{eqnarray}} 
\newcommand{\eqarraystar} {\begin{eqnarray*}} 
\newcommand{\enarraystar} {\end{eqnarray*}} 
\newcommand{\bel}{\begin{lemma}} 
\newcommand{\el}{\end{lemma}} 
\newcommand{\bes}{\begin{step}} 
\newcommand{\es}{\end{step}} 
\newcommand{\bpr}{\begin{proof}} 
\newcommand{\epr}{\end{proof}}

\renewcommand{\section}{\secdef\sct\sect} 
\newcommand{\sct}[2][default]{\refstepcounter{section} 
\vspace{0.8cm} 
\setcounter{equation}{0} 
\centerline{ 
\large\scshape \arabic{section}.\ #1} 
\vspace{0.2cm}} 
\newcommand{\sect}[1]{ 
\vspace{0.8cm} 
\centerline{\large\scshape #1} 
\vspace{0.2cm}} 
 
\renewcommand{\subsection}{\secdef \subsct\sbsect} 
\newcommand{\subsct}[2][default]{\refstepcounter{subsection} 
\nopagebreak 
\vspace{0.5\baselineskip} 
{\flushleft\bf \arabic{section}.\arabic{subsection}~\bf #1  } 
\nopagebreak} 
\newcommand{\sbsect}[1]{\vspace{0.1cm}\noindent 
{\bf #1}\vspace{0.1cm}}

\renewcommand{\subsubsection}{%
\secdef \subsubsect\sbsbsect} 
\newcommand{\subsubsect}[2][default]{%
\refstepcounter{subsubsection} 
\nopagebreak 
\vspace{0.1\baselineskip} 
\nopagebreak 
{\flushleft 
\sffamily\slshape 
\arabic{section}.\arabic{subsection}.\arabic{subsubsection} 
\ %
\sffamily #1\/.}\ } 
\newcommand{\sbsbsect}[1]{\vspace{0.1cm}\noindent 
{\bf #1}\ } 
 
\renewcommand{\d}{{\rm d}}

\newcommand{\eps}{\varepsilon}

\newcommand{\supp}{{\operatorname {supp}}} 
\newcommand{\dist}{{\operatorname {dist}}}

\newcommand{\Ccal}   {{\mathcal C }}

\newcommand{\Mcal}   {{\mathcal M }} 
\newcommand{\Ncal}   {{\mathcal N }}

\newcommand{\Xcal}   {{\mathcal X }}

\begin{document} 
 
\title[Large systems of path-repellent Brownian motions]{\large 
Large systems of path-repellent Brownian motions \medskip in a trap at positive temperature} 
 
\author[Stefan Adams, Jean-Bernard Bru and Wolfgang 
        K{\"o}nig]{} 
\maketitle

\thispagestyle{empty} 
\vspace{0.2cm} 
 
\centerline {\sc By Stefan Adams\footnote{Max-Planck Institute for Mathematics in the Sciences, Inselstra{\ss}e 22-26, D-04103 Leipzig, Germany, {\tt adams@mis.mpg.de}}, Jean-Bernard Bru\footnote{Fachbereich Mathematik und Informatik,
Johannes-Gutenberg-Universit\"at Mainz, Staudingerweg 9, D-55099 Mainz, Germany, {\tt  jbbru@mathematik.uni-mainz.de}} \/ and  Wolfgang
K{\"o}nig\footnote{Mathematisches Institut, Universit\"at Leipzig, Augustusplatz 10/11, D-04109 Leipzig, Germany, {\tt koenig@math.uni-leipzig.de}}}
\vspace{0.4cm}
\renewcommand{\thefootnote}{}
 
\vspace{.5cm} 
 
\centerline{\small(\version)} 
\vspace{.5cm} 
 

\bigskip

\begin{quote} 
{\small {\bf Abstract:}} We study a model of $ N $ mutually repellent Brownian motions under confinement to stay in some bounded region of space. Our model is defined in terms of a transformed path measure under a trap Hamiltonian, which prevents the motions from escaping to infinity, and a pair-interaction Hamiltonian, which imposes a repellency of the $N$ paths. In fact, this interaction is an $N$-dependent regularisation of the Brownian intersection local times, an object which is of independent interest in the theory of stochastic processes. 

The time horizon (interpreted as the inverse temperature) is kept fixed. We analyse the model for diverging number of Brownian motions in terms of a large deviation principle. The resulting variational formula is the positive-temperature analogue of the well-known Gross-Pitaevskii formula, which approximates the ground state of a certain dilute large quantum system; the kinetic energy term of that formula is replaced by a probabilistic energy functional. 

This study is a continuation of the analysis in \cite{ABK04} where we considered the limit of diverging time (i.e., the zero-temperature limit) with fixed number of Brownian motions, followed by the limit for diverging number of motions.
\end{quote}

\vfill


 
\eject 
 
\setcounter{section}{0} 
\section{Introduction and main results}\label{Intro}

\subsection{Background.}

\noindent One of the important problems in mathematical physics is the understanding of large systems of interacting quantum particles at extremely low or zero temperature. This question is raised in the literature since decades and is studied under a lot of different view points. 
Large systems of {\it bosons\/} (i.e., quantum particle systems whose wave functions are invariant under permutation of the single-particle variables), at extremely low temperatures, may undergo a phase transition, the so-called {\it Bose-Einstein condensation}: a macroscopic portion of the system is described by one suitable single-particle wave function. This phenomenon was shown for ideal (i.e., non-interacting) gases already in 1925 by Bose and Einstein. A rigorous understanding of condensation effects in various ultra cold materials remains a widely open and challenging problem until today. This question received an enormous impetus by the first experimental realisation of Bose-Einstein condensation in 1995. 

A many-particle quantum system is mathematically described by a $N$-particle {\it Hamilton} operator $H_N$ consisting of a kinetic energy term, a trap term and an interaction term. Its spectral analysis, at least for realistic interacting models, is out of reach of contemporary analysis. Rigorous theoretical research started with Bogoliubov and Landau in the 1940ies, followed by Penrose, Feynman and many others. They analysed simplified mathematical models featuring only the most important physical phenomena. However, these approaches turned out to be intuitively appealing and relevant. See \cite{AB04a,AB04b} for a review and some recent results. 

Another mathematical approach is to consider systems that are dilute on a particular scale and are kept within a bounded region by the presence of a trap. Here \lq dilute\rq\ means that the range of the interparticle interaction is small compared with the mean particle distance. These systems are supposed to be easier to analyse at least as it concerns the ground state. In a particular dilute situation, the ground states and their energy were analysed in the many-particle limit \cite{LSSY05}. It turned out that the well-known {\it Gross-Pitaevskii formula\/} describes the system remarkably well. This variational formula has a kinetic term (the usual energy), a trap term and a quartic term with a pre factor. As was predicted by earlier theoretical work, the only parameter of the pair interaction functional that persists in the limit is its scattering length. See \cite{PS03} for an overview about the physics and \cite{LSSY05} for an account on recent mathematical research.

However, the mathematically rigorous understanding of large quantum systems at {\it positive\/} temperature is still incomplete. For dilute systems of {\it fermions\/} (i.e., quantum particle systems whose wave functions are antisymmetric under permutation of the single-particle variables), first results for positive temperature are in \cite{Sei05}. One main concern of quantum statistical mechanics is to evaluate the trace of the Boltzmann factor $ {\rm e}^{-\beta H_N} $ for inverse temperature $ \beta > 0 $ to calculate all thermodynamic functions. The {\it Feynman-Kac formula\/} provides a representation of these traces as functional integrals over the space of Brownian paths on the finite time horizon $[0,\beta]$ \cite{Gin71}. Hence, it is clear that an appropriate description of quantum systems at positive temperature is given in terms of independent Brownian motions in a trap with a mutually repellent pair interaction. The trap and the interaction are imposed via exponential densities, so-called Hamiltonians. 

In the present paper, we make a contribution to a rigorous analysis of a certain model of a large number of mutually repellent Brownian motions in a trap at any positive temperature. We introduced this model in earlier work \cite{ABK04}. The pair interaction in that model is a {\it path\/} interaction, not a {\it particle\/} interaction. It turned out there that its behaviour in the zero-temperature limit is asymptotically well described by a variational formula known as the {\it Hartree formula}. Therefore, we call this Brownian model the {\it Hartree model}. The interaction Hamiltonian is given via a double time integration and thus the Hartree model is related to Polaron type models \cite{DV83}, \cite{BDS93}, where instead of several paths a single path is considered. In \cite{ABK04} we showed that the many-particle limit of the Hartree formula is well approximated by the above mentioned Gross-Pitaevskii formula. However, the decisive parameter here is not the scattering length, but the {\it integral\/} of the pair interaction functional. 

While that result describes the zero-temperature situation, in the present work we study the case of {\it positive\/} temperature, i.e., interacting Brownian motions on a fixed {\it finite\/} time horizon. Our main result is a description of the many-particle limit in terms of a variational formula that is analogous to the Gross-Pitaevskii formula. The only difference is the kinetic energy term, which is now replaced by a {\it probabilistic energy\/} term. In the language of large deviation theory, this term is the rate function that governs a certain large deviation principle. In our particular case it is a certain Legendre-Fenchel transform of an exponential Brownian integral.

The system of motions is dilute on the same scale as the above mentioned models. This means in our particular case that the repellent interaction is an approximation of a highly irregular functional of the motions, the so-called {\it Brownian intersection local times}, which measure the amount of time that is spent by two motions at their intersection points. This object is of independent interest in the theory of stochastic processes. Our main result implicitly states a large deviation principle for the mean of regularisations of the intersection local times, taken over all mutual intersections of a large number of Brownian paths.

The remainder of Section~\ref{Intro} is organised as follows. We introduce the model in Section~\ref{sec-model} and present our main result and some conclusions in Section~\ref{sec-result}. In Section~\ref{sec-RevI} we embed these results in a broader perspective, discuss our results and mention some open problems.

The remainder of the paper is structured as follows. In Section~\ref{sec-Vari} we prove some properties of the probabilistic energy term and the variational formula. Section~\ref{sec-ProofN} contains the proof of our main result. In the Appendix we give a short account on large deviation theory in Section~\ref{LDP} and recall a related result by Lieb {\it et~al.} on the large-$N$ limit of the ground state in Section~\ref{GPtheory}.

\subsection{The model.}\label{sec-model} 

\noindent We consider a family of  $N$ independent Brownian motions, $(B_t^{\ssup{1}})_{t\geq 0}, \dots,(B_t^{\ssup{N}})_{t\geq 0}$, in $\R^d$ with generator $-\Delta$ each. We assume that each motion possesses the same initial distribution, which we do not want to make explicit. The model we study  is defined in terms of a Hamiltonian which consists of two parts: a trap part, 
\begin{equation}\label{ham1}
H_{N,\beta}=\sum_{i=1}^N\int_0^\beta W(B_s^{\ssup{i}})\,\d s, 
\end{equation}
and a pair-interaction part, 
\begin{equation}\label{ham3}
K_{N,\beta}=\sum_{1\leq i<j\leq N}\frac 1\beta\int_0^\beta\int_0^\beta v\bigl(|B_s^{\ssup{i}}-B_t^{\ssup{j}}|\bigr)\,\d s\d t.
\end{equation}
Here $W\colon\R^d\to[0,\infty]$ is the so-called {\it trap potential} satisfying $\lim_{|x|\to\infty}W(x)=\infty$, and $v\colon (0,\infty)\to[0,\infty]$ is a 
pair-interaction function satisfying $0<\lim_{r\downarrow 0}v(r)\leq\infty$ and $\int_{\R^d}v(|x|)\,\d x<\infty$. We are interested in the large-$N$ behaviour of the transformed path measure,
\begin{equation}\label{ZNKdef} 
{\rm e}^{-H_{N,\beta}-K_{N,\beta}}\,\d\P.
\end{equation} 
Here $\beta\in(0,\infty)$ is a finite time horizon which we will keep fixed in this paper. The trap part effectively keeps the motions in a bounded region of the space $\R^d$. 
Through the pair interaction $K_{N,\beta}$, the $i$-th Brownian motion interacts with the mean of the whole path of the $j$-th motion, taken over all times before $\beta$. Hence,  the interaction is not a {\it particle\/} interaction, but a {\it path\/} interaction. We are most interested in the case $\lim_{r\downarrow 0}v(r)=\infty$, where the pair-interaction repels all the motions from each other (more precisely, their paths). In order to keep the notation simpler, we abstained from normalising the path measure in \eqref{ZNKdef}. 

The model in \eqref{ZNKdef} was introduced and studied in \cite{ABK04}; see Section~\ref{sec-RevI} for results from that paper and a discussion of the physical relevance of the model. In particular, in the limit $\beta\to\infty$, followed by $N\to\infty$, a certain variational formula appears that is called the {\it Hartree formula\/} in the literature. Therefore, we call the model in \eqref{ZNKdef} the {\it Hartree model}.

When we take the limit as $N\to\infty$, we will not keep the pair-interaction function $v$ fixed, but we replace it by the rescaled version $v_N(\cdot)=N^{d-1}v(N\,\cdot)$. In other words, we replace $K_{N,\beta}$ with
\begin{equation}\label{KNbetaNdef}
K_{N,\beta}^{\ssup N}=\frac 1N\sum_{1\leq i<j\leq N}\frac 1\beta\int_0^\beta\int_0^\beta N^{d}v\bigl(N|B_s^{\ssup{i}}-B_t^{\ssup{j}}|\bigr)\,\d s\d t.
\end{equation}
Note that $N^{d}v(N\,\cdot)$ is an approximation of the Dirac measure at zero times the integral of $v\circ|\cdot|$, hence the double integral in \eqref{KNbetaNdef} is an approximation of the {\it Brownian intersection local times\/} at zero, an important object in the present paper. The Brownian intersection local times measure the time spent by two motions on the intersection of the their paths, see Section~\ref{sec-BISLT}. A natural sense can be given to this object only in dimensions $d\in\{2,3\}$. Therefore, our analysis is naturally restricted to these dimensions. The main difficulty in the analysis of the model will stem from the $N$-dependence and the high irregularity of the pair-interaction part.

\subsection{Results.}\label{sec-result}

\noindent We now formulate our results on the behaviour of the Hartree model in the limit as $N\to\infty$, with $\beta>0$ fixed. First we introduce an important functional, which will play the role of a probabilistic energy functional. Define $J_\beta\colon \Mcal_1(\R^d)\to[0,\infty]$ as the Legendre-Fenchel transform of the map $\Ccal_{\rm b}(\R^d)\ni f\mapsto \frac{1}{\beta}\log \E[{\rm e}^{\int_0^\beta f(B_s)\,\d s}]$ on the set $\Ccal_{\rm b}(\R^d)$ of continuous bounded functions on $\R^d$, where $ (B_s)_{s\ge 0} $ is one of the above Brownian motions. That is,
\begin{equation}\label{JTKdef}
J_\beta(\mu)=\sup_{f\in \Ccal_{\rm b}(\R^d)}\Bigl(\langle\mu,f\rangle-\frac{1}{\beta}\log \E\big[{\rm e}^{\int_0^\beta f(B_s)\,\d s}\big]\Bigr),\qquad\mu\in\Mcal_1(\R^d).
\end{equation}
Here $\Mcal_1(\R^d)$ denotes the set of probability measures on $\R^d$. Note that $ J_\beta $ depends on the initial distribution of the Brownian motion. In Lemma~\ref{Jnotinfty} below we show that $ J_{\beta} $ is not identical to $ +\infty $. Alternate expressions for $J_\beta$ are given in Lemma~\ref{Jbetaident} below. Clearly, $J_\beta$ is a lower semi continuous and convex functional on $\Mcal_1(\R^d)$, which we endow with the topology of weak convergence induced by test integrals against continuous bounded functions. However, $J_\beta $ is {\it not\/} a quadratic form coming from any linear operator. We wrote $\langle \mu,f\rangle=\int_{\R^d}f(x)\,\mu(\d x)$ and use also the notation $\langle f,g\rangle=\int_{\R^d}f(x)g(x)\,\d x$ for integrable functions $f,g$. If $\mu$ possesses a Lebesgue density $\phi^2$ for some $L^2$-normalised $\phi\in  L^2$, then we also write $J_\beta(\phi^2)$ instead of $J_\beta(\mu)$. In Lemma~\ref{densitiesofmeasuresJ} below it turns out that $J_\beta(\mu)=\infty$ if $\mu$ fails to have a Lebesgue density.

In the language of the theory of large deviations, $J_\beta$ is the rate function that governs a  certain large deviation principle. (See Section~\ref{LDP} for the notion and some remarks on large deviation theory.) The object that satisfies this principle is the mean of the $N$  normalised occupation measures, 
\begin{equation}\label{meanmu}
\overline\mu_{N,\beta}=\frac 1N\sum_{i=1}^N \mu_\beta^{\ssup i},\qquad N\in\N.
\end{equation}
Here
\begin{equation}\label{muTidef}
\mu_\beta^{\ssup{i}}(\d x)=\frac 1\beta\int_0^\beta\delta_{B_s^{\ssup{i}}}(\d x)\,\d s,\qquad i=1,\dots,N,
\end{equation}
is the {\it normalised occupation measure\/} of the $i$-th motion, which is a random element of $\Mcal_1(\R^{d})$. It measures the time spent by the $i$-th Brownian motion  in a given region. One can write the Hamiltonians in terms of the occupation measure as
\begin{equation}
H_{N,\beta}=\beta\sum_{i=1}^N\langle W,\mu_\beta^{\ssup{i}}\rangle\qquad\mbox{and}\qquad K_{N,\beta}=\beta\sum_{1\leq i<j\le N}\langle \mu_\beta^{\ssup{i}}, {V}\mu_\beta^{\ssup{j}}\rangle,
\end{equation}
where we denote by $V$ the integral operator with kernel $v$, i.e., ${V}f(x)=\int_{\R^d}v(|x-y|)f(y)\,\d y$ and analogously for measures $V\mu(x)=\int_{\R^d}\mu(\d y)\,v(|x-y|)$. Hence, it is natural to expect that the asymptotic of the Hartree model can be expressed in terms of asymptotic properties of $\overline\mu_{N,\beta}$. 

We are heading towards a formulation of our main result. Our precise assumptions on the trap potential, $W$, and on the pair-interaction functional, $v$, are the following.

\medskip

\noindent{\bf Assumption (W).} {\it $W\colon\R^d\to[0,\infty]$ is continuous in $\{W<\infty\}$ with $\lim_{R\to\infty}\inf_{|x|>R}W(x)=\infty$. Furthermore, $\{W<\infty\}$ is either equal to $\R^d$ or is a bounded connected open set.}

\medskip

\noindent{\bf Assumption (v).} {\it $v\colon[0,\infty)\to[0,\infty]$ is measurable, $\int_{\R^d}v(|x|)\,\d x<\infty$ and $\int_{\R^d}v(|x|)^2\,\d x<\infty$.}

\medskip

In order to avoid trivialities, we tacitly assume that the support of the initial distribution of the Brownian motions is contained in the set $\{W<\infty\}$.

Now we formulate our main result. As we already indicated, the main role in the analysis of the Hartree model is played by the mean of the normalised occupation measures in \eqref{meanmu}.

\begin{theorem}[Many-particle limit for the Hartree model]\label{Kmodel,N} Assume that $d\in \{2,3\}$ and let $W$ and $v$ satisfy Assumptions~(W) and (v), respectively. Introduce
\begin{equation}\label{vnorm}
 \sca(v):=\frac 1{8\pi}\int_{\R^d} v(|y|)\,\d y <\infty.
\end{equation}  
Fix $\beta>0$. Then, as $N\to\infty$, the mean $ \overline{\mu}_{N,\beta}=\frac{1}{N}\sum_{i=1}^N\mu_\beta^{\ssup{i}} $ of the normalised occupation measures  satisfies a large deviation principle on $\Mcal_1(\R^d)$ under the measure with density ${\rm e}^{-H_{N,\beta}-K_{N,\beta}^{\ssup N}}$ with speed $ N \beta$ and rate function 
\begin{equation}\label{Ibetadef}
I^{\ssup{\otimes}}_\beta(\mu)=
\begin{cases} J_\beta(\phi^2)+\langle W,\phi^2\rangle+4\pi{\sca}(v)\,||\phi||_4^4&\text{if }\phi^2=\frac {\d\mu}{\d x}\text{ exists,}\\
\infty&\text{otherwise.}
\end{cases}
\end{equation}
The level sets $ \{\mu\in\Mcal_1(\R^d)\colon I^{\ssup{\otimes}}_\beta(\mu)\le c\}$, $c\in\R$, are compact.
\end{theorem}

We also write $I^{\ssup{\otimes}}_\beta(\phi^2)$ if $\phi^2=\frac {\d\mu}{\d x}$. To be more explicit, the large deviation principle for $ \overline{\mu}_{N,\beta}$ means that
\begin{equation}\label{LDPmuN}
\lim_{N\to\infty}\frac 1{N\beta}\log\E\big[{\rm e}^{-H_{N,\beta}-K_{N,\beta}^{\ssup N}}\1_{\overline{\mu}_{N,\beta}\in \,\bdot\,}\big]=-\inf_{\phi^2\in\, \bdot\,}I^{\ssup{\otimes}}(\phi^2)\qquad\mbox{weakly},
\end{equation}
where we identify $\Mcal_1(\R^d)$ with the unit sphere in $L^2(\R^d)$ via the relation $\phi^2(x)\,\d x=\mu(\d x)$. The convergence in \eqref{LDPmuN} is in the weak sense, i.e., the lower bound holds for open sets and the upper bound for closed sets (see Section~\ref{LDP} for more details). Here we refer to the weak topology on $ \Mcal_1(\R^d) $. See Section~\ref{sec-BISLT} for a heuristic explanation of the assertion of Theorem~\ref{Kmodel,N}.

In Assumption (v) we require that $ v\circ|\cdot|\in L^2(\R^d) $. This is  needed in our proof of the {\it lower\/} bound in (\ref{LDPmuN}) only. We think that this assumption is technical only and could be relaxed if higher integrability properties of the elements of the level sets of $ I^{\ssup{\otimes}}_\beta $ were known.

Interesting conclusions of Theorem~\ref{Kmodel,N} are as follows. 
For $\sca>0$, we introduce the variational formula
\begin{equation}\label{GPTK}
\chi_{\sca}^{\ssup{\otimes}}(\beta)=\inf_{\phi\in L^2(\R^d)\colon \|\phi\|_2=1}\Big( J_\beta(\phi^2)+\langle W,\phi^2\rangle +4\pi\sca\,||\phi||_4^4\Big),
\end{equation}
which, for $ \sca=\sca(v) $, is the minimum of the rate function $I_\beta^{\ssup{\otimes}}$ defined in \eqref{Ibetadef}. Here are some facts about the minimiser in \eqref{GPTK}.

\begin{lemma}[Analysis of $\chi_{\sca}^{\ssup{\otimes}}(\beta)$]\label{GPTKmin}
Fix $\beta> 0 $ and $\sca>0$. 
\begin{enumerate}
\item[(i)] There exists a unique $L^2$-normalised minimiser $\phi_*\in L^2(\R^d)\cap L^4(\R^d) $ of the right hand side of \eqref{GPTK}.

\item[(ii)] For any neighbourhood $\Ncal\subset L^2(\R^d)\cap L^4(\R^d)$ of $\phi_*$,
$$
\inf_{\phi\in L^2(\R^d)\colon \|\phi\|_2=1,\phi\notin\Ncal} \Bigl( J_\beta(\phi^2)+\langle W,\phi^2\rangle +4\pi\sca||\phi||_4^4\Bigr) >\chi_{\sca}^{\ssup{\otimes}}(\beta).
$$
Here \lq neighbourhood\rq\ refers to any of the three following topologies: weakly in $L^2$, weakly in $L^4$, and weakly in the sense of probability measures, if $\phi$ is identified with the measure $\phi(x)^2\,\d x$. 
\end{enumerate}
\end{lemma}

Now we can state some conclusions about the large-$N$ behaviour of the total expectation of the exponential Hamiltonian and about a kind of Law of Large Numbers. The proof is simple and omitted.
 
\begin{cor}
Let the assumptions of Theorem~\ref{Kmodel,N} be satisfied. Then the following holds.
\begin{enumerate}
\item[(i)] 
\begin{equation}\label{meanfreeenergy} 
\lim_{N\to\infty}\frac 1{\beta N}\log\E\big[{\rm e}^{-H_{N,\beta}-K^{\ssup N}_{N,\beta}}\big]=-\chi^{\ssup{\otimes}}_{{\sca}(v)}(\beta). 
\end{equation} 
\item[(ii)] As $ N\to\infty $,  $ \overline{\mu}_{N,\beta}$ converges in distribution under the measure with density ${\rm e}^{-H_{N,\beta}-K_{N,\beta}^{\ssup N}}$ towards the measure $\phi_*(x)^2\,\d x$, where $\phi_*\in L^2(\R^d)$ is the unique minimiser in \eqref{GPTK} with $\alpha=\alpha(v)$ as defined in \eqref{vnorm}.
\end{enumerate}
\end{cor}

\subsection{Relation with quantum statistical mechanics.}\label{sec-RevI}

\noindent In this section we briefly explain the relation between the Hartree model in \eqref{ZNKdef} and quantum statistical mechanics.
  
An $N$-particle quantum system is described by the $N$-particle Hamilton operator
\begin{equation}\label{HNdef}
H_N=\sum_{i=1}^N\Bigl[-\Delta_i^2 +W(x_i)\Bigr]+\sum_{1\le i< j\le N}v(|x_i-x_j|),\quad x=(x_1,\ldots,x_N)\in\R^{dN}.
\end{equation}
For $\beta>0$ fixed, the trace of the Boltzmann factor, ${\rm e}^{-\beta H_N}$, is given, via the Feynman-Kac formula, by a Brownian model similar to the one in \eqref{ZNKdef}, where the Brownian motions are conditioned to terminate at their starting points (Brownian bridges) and the initial measure is the Lebesgue measure \cite{Gin71}, \cite{BR97}. However, the interaction Hamiltonian is, instead of $K_{N,\beta}$ in \eqref{ham3},
\begin{equation}\label{canonical}
\sum_{1\le i<j\le N}\int_0^{\beta} v\big(|B_s^{\ssup{i}}-B_s^{\ssup{j}}|\big)\,\d s.
\end{equation}
This is a {\it particle\/} interaction involving only one time axis for all the motions, in contrast to the time-pair integration in \eqref{ham3}. Note that there is no Hamilton operator such that the total mass of the Hartree model is equal to the trace of the corresponding Boltzmann factor. The Hartree model features the mutually repellent nature of the trace of ${\rm e}^{-\beta H_N}$ in a form which is more accessible to a rigorous stochastic analysis. The study of the large-$N$ behaviour of that trace, i.e., of the Brownian model with the interaction in \eqref{canonical}, is deferred to future work.

For describing large systems of {\it bosons\/} at positive temperature, one has to consider the trace of the projection of $H_N$ to the subspace of symmetric wave functions. The corresponding Brownian model is given in terms of Brownian bridges with symmetrised initial and terminal conditions. The effect of this symmetrisation on the large-$N$ limit is studied elsewhere \cite{AK06}.

Let us now comment on the physical relevance of the Hartree model. Its usefulness as a simplified model for the ground state of $ H_N $ is well-known, see the physics monograph \cite[Ch.~12]{DN05}. The precise relation is as follows. Recall the integral operator, $V$, with kernel $v$, and introduce the variational formula
\begin{equation}\label{chivKdef}
\chi_{N }^{\ssup{\otimes}}=\frac 1N\inf_{h_1,\dots,h_N\in H^1(\R^d)\colon \|h_i\|_2=1\,\forall i}\Big(\sum_{i=1}^N\Bigl[\|\nabla h_i\|_2^2+\langle W,h_i^2\rangle\Bigr]+\sum_{1\leq i<j\leq N} \langle h_i^2, Vh_j^2\rangle\Big).
\end{equation}
In \cite{ABK04} we showed that there exist tuples $(h_1,\dots,h_N)$ of minimisers, and we characterised them in terms of Euler-Lagrange equations and proved some regularity properties. Note that there is no convexity argument available, which leads us to the conjecture that the tuples of minimisers are not unique. Note also that 
\begin{equation}\label{chivKident}
\chi_{N }^{\ssup{\otimes}}=\frac 1N\inf_{h_1,\dots,h_N\in H^1(\R^d)\colon \|h_i\|_2=1\,\forall i}\langle h,H_N h\rangle,\qquad  \mbox{where }h=h_1\otimes\cdots\otimes h_N.
\end{equation}
Hence, one can conceive $\chi_{N }^{\ssup{\otimes}}$ as the {\it ground product-state energy} of $ H_N$, i.e., as the ground state energy of the restriction of $ H_N$ to the set of $N$-fold product states. If $(h^{\ssup{N}}_1,\dots,h^{\ssup{N}}_N)$ is a minimiser, we call  $h^{\ssup{N}}:=h_1^{\ssup{N}}\otimes\cdots\otimes h_N^{\ssup{N}} $ a {\it ground product-state}. One of the main results of \cite{ABK04}, Th.~1.7, states that, for any fixed $N\in\N$,
\begin{equation}\label{betatoinfty}
\lim_{\beta\to\infty}\frac1{N\beta}\log\E\big[{\rm e}^{-H_{N,\beta}-K_{N,\beta}}\big]=-\chi_{N }^{\ssup{\otimes}}.
\end{equation}
The proof shows that the tuple of normalised occupation measures, $(\mu_\beta^{\ssup{1}},\dots, \mu_\beta^{\ssup{N}})$ (recall \eqref{muTidef}) stands in a one-to-one relation with the minimiser tuples $(h_1,\dots,h_N)$ of \eqref{chivKdef}, in the sense of a large deviation principle, analogously to \eqref{LDPmuN}. This result illustrates the close connection between the zero-temperature Hartree model and the ground states of the Hamilton operator $H_N$. One main novelty of the present paper is to demonstrate that also at positive temperature the Hartree model is a useful simplification of the trace of $ {\rm e}^{-\beta H_N} $.

It is instructive to compare the main result of the present paper, Theorem~\ref{Kmodel,N}, to the zero-temperature analogue of that result, which we derived in \cite{ABK04} and which served us as a main motivation for the present work. It turned out there that the {\it Gross-Pitaevskii formula\/} well approximates the ground product-state energy $\chi_N^{\ssup{\otimes}}$ in the large-$N$ limit, provided that the interaction functional $v$ is rescaled in the same manner as in \eqref{KNbetaNdef}. The Gross-Pitaevskii formula, derived independently by Gross and Pitaevskii in 1961 on the basis of the method initiated by Bogoliubov and Landau in the 1940ies, has a parameter $\sca>0$ and is defined as follows.
\begin{equation}\label{GP}
\chi_\sca^{\ssup{\rm GP}}=\inf_{\phi\in H^1(\R^d)\colon\|\phi\|_2=1}\Big(\|\nabla\phi\|_2^2+\langle W,\phi^2\rangle+4\pi \sca \|\phi\|_4^4\Big).
\end{equation}
There is a unique minimiser $\phi^{\ssup{\rm GP}}_{\sca}$, which is positive and continuously differentiable with H\"older continuous derivatives  \cite{LSSY05}.

The main result of the present paper, Theorem~\ref{Kmodel,N}, see \eqref{betatoinfty}, is the positive-temperature analogue of the following result, which is \cite[Th.~1.14]{ABK04}.

\begin{theorem}[Large-$N$ asymptotic of $\chi^{\ssup{\otimes}}_{N }$]\label{thm-chiNconv} Let $d\in\{2,3\}$. Assume that $v$ satisfies Assumption (v). Replace $v$ by $v_N(\cdot)=N^{d-1}v(N\,\cdot)$. Let $(h_1^{\ssup{N}},\dots,h_N^{\ssup{N}})$ be any minimiser on the right hand side of \eqref{chivKdef} (with $v$ replaced by $v_N(\cdot)=N^{d-1}v(N\,\cdot)$).
Define $\phi_N^2=\frac 1N\sum_{i=1}^N (h_i^{\ssup{N}})^2$. Then we have
$$
\lim_{N\to\infty}\chi^{\ssup{\otimes}}_{N }=\chi^{\ssup{\rm GP}}_{ \sca(v)}\qquad\mbox{and}\qquad \phi_N^2\to\big(\phi^{\ssup{{\rm GP}}}_{ \sca(v)}\big)^2,
$$
where $\sca(v)$ is the integral introduced in \eqref{vnorm}. The convergence of $\phi_N^2$ is in the weak $L^1(\R^d)$-sense and weakly for the probability
measures $\phi_n^2(x)\,\d x$ towards the measure $(\phi^{\ssup{{\rm GP}}}_{ \sca(v)})^2(x)\,\d x$.
\end{theorem}

We remark that Lieb {\it et~al.} proved the analogous result for the ground state of $H_N$, see Section~\ref{GPtheory} in the appendix. Note that in $d=2$ the scaling of $v$ differs from the one used in Theorem~\ref{thm-chiNconv}. Moreover, the parameter $\sca(v)$ in Theorem~\ref{thm-chiNconv} is replaced by the scattering length of $v$ in the result of Lieb {\it et~al.} The integral $\sca(v)$ is known as the {\it first Born approximation\/} of the scattering length of $v$ \cite{LSSY05}.

We conjecture that $\chi_{\sca}^{\ssup{\otimes}}(\beta)$ in \eqref{GPTK} converges to the Gross-Pitaevskii formula as $\beta\to\infty$. A proof of this is due to future work.

\section{Variational analysis}\label{sec-Vari}

In this section we derive some useful properties of the probabilistic energy functional $J_\beta $ introduced in \eqref{JTKdef} in Section~\ref{sec-Jbeta}, and prove the existence and uniqueness of minimisers in the variational formula $\chi^{\ssup{\otimes}}_\sca(\beta)$ introduced in \eqref{GPTK}  in Section~\ref{sec-chi}.

\subsection{Some properties of $\boldsymbol{J_\beta}$.}\label{sec-Jbeta}

\noindent  First we show that $J_\beta$ is not identically equal to $+\infty$.

\begin{lemma}\label{Jnotinfty}
There is $\mu\in\Mcal_1(\R^d)$ such that $J_\beta(\mu)<\infty$.
\end{lemma}

\begin{proofsect}{Proof} Recall that  $J_\beta$ is the Legendre-Fenchel transform of the map $\Ccal_{\rm b}(\R^d)\ni f\mapsto\log \E[{\rm e}^{\beta\langle f,\mu_{\beta}\rangle}]$, where we recall that $\mu_\beta$ is the normalised occupation measure of one of the Brownian motions. Recall the mean of the $N$ normalised occupation measures from \eqref{meanmu}. Now pick a continuous function $g\colon\R^d\to[0,\infty)$ satisfying $\lim_{R\to\infty}\inf_{|x|\geq R}g(x)=\infty$. Then we have, for any $C>0$, by splitting the probability space into $\{\langle g,\overline\mu_{N,\beta}\rangle\leq C\}$ and its complement,
$$
\begin{aligned}
-\infty&<\log\E\big[{\rm e}^{-\langle g,\mu_{\beta}\rangle}\big]=\limsup_{N\to\infty}\frac1N\log\E\big[{\rm e}^{-N\langle g,\overline\mu_{N,\beta}\rangle}\big]\\
&\leq \max\Big\{-C,\limsup_{N\to\infty}\frac1N\log\E\big[{\rm e}^{-N\langle g,\overline\mu_{N,\beta}\rangle}\1_{\{\langle g,\overline\mu_{N,\beta}\rangle\leq C\}}\big]\Big\}.
\end{aligned}
$$
According to \cite[Th.~4.5.3(b)]{DZ98}, the sequence $(\overline\mu_{N,\beta})_{N\in\N}$ satisfies the upper bound in the large deviation principle for compact sets with rate function equal to $J_\beta$. By Prohorov's Theorem, the set $\{\mu\in\Mcal_1(\R^d)\colon \langle g,\mu\rangle\leq C\}$ is compact. Furthermore, note that the map $\mu\mapsto -\langle g,\mu\rangle$ is upper semi continuous.
Hence, the upper-bound part in Varadhan's Lemma, \cite[Lemma 4.3.6]{DZ98}, implies that
$$
\limsup_{N\to\infty}\frac1N\log\E\big[{\rm e}^{-N\langle g,\overline\mu_{N,\beta}\rangle}\1_{\{\langle g,\overline\mu_{N,\beta}\rangle\leq C\}}\big]
\leq -\inf_{\mu\in\Mcal_1(\R^d)\colon \langle g,\mu\rangle\leq C}\Big(\langle g,\mu\rangle +J_\beta(\mu)\Big).
$$
Picking $C$ large enough, we find that
$$
\infty>\inf_{\mu\in\Mcal_1(\R^d)\colon \langle g,\mu\rangle\leq C}\Big(\langle g,\mu\rangle +J_\beta(\mu)\Big).
$$
This implies that $J_\beta$ is not identically equal to $\infty$.
\qed
\end{proofsect}

Now we show that $J_\beta$ is infinite in any probability measure in $ \R^d$ that fails to have a Lebesgue density.

\begin{lemma}\label{densitiesofmeasuresJ}
If $\mu\in\Mcal_1(\R^d)$ is not absolutely continuous, then $J_\beta(\mu)=\infty$.
\end{lemma}
\begin{proofsect}{Proof}
We write $\lambda$ for the Lebesgue measure on $\R^d$. Pick $\mu\in\Mcal_1(\R^d)$ that is not absolutely continuous. Then there is a Borel set $ A\subset\R^d $ such that $ \lambda(A)=0 $ and $ \mu(A)>0 $. Let $ M>0 $. We show that $ J_\beta(\mu) \ge M $. Pick $ K=\frac{4M}{\mu(A)} $ and $ \eta=\frac{1}{2}\mu(A) $. We may assume that $ 2\le {\rm e}^{\beta\frac{K}{4}\mu(A)} $. Let $(Q_{\eps})_{\eps >0}$ be  an increasing family of open subsets of $\R^d $ such that $ A\subset Q_{\eps} $ and $ \lambda(Q_{\eps})<\eps $ for any $\eps>0$. Let $\mu_\beta$ is the (random) normalised occupation measure of a Brownian motion $ (B_s)_{s\ge 0}$. Pick $ \eps>0 $ and $ \delta>0 $ so small that 
$$
 \P\Bigl(\mu_\beta(U_{\delta}(Q_{\eps}))>\eta\Bigr) <{\rm e}^{-\beta K(1-\eta)}, 
$$  where $ U_{\delta}(Q_{\eps}) $ is the $ \delta$-neighbourhood of $ Q_{\eps} $.
This is possible since
$$
\limsup_{\eps\downarrow 0}\limsup_{\delta\downarrow 0}\mu_\beta(U_{\delta}(Q_{\eps}))=\limsup_{\eps\downarrow 0}\mu_\beta(\overline{Q}_{\eps})=\mu_\beta\Bigl(\bigcap_{\eps >0}
\overline{Q}_{\eps}\Bigr)=0\, \mbox{ a.s.}
$$
Now we pick a function $ f\in\Ccal_{\rm b}(\R^d) $ such that $ 0\le f\le K$, $\mbox{ supp }(f)\subset U_{\delta}(Q_{\eps}) $ and $ f|_{Q_{\eps}}=K $. Then $ \langle f,\mu\rangle \ge \int_{Q_{\eps}} f(x)\,\mu(\d x)\ge K\mu(A)$. Furthermore, 
\begin{eqnarray}
\begin{aligned}
\E\Big[{\rm e}^{\int_0^\beta f(B_s)\,\d s}\Big]
&= \E\Big[{\rm e}^{\beta\langle f,\mu_\beta\rangle}\1\{\mu_\beta(U_{\delta}(Q_{\eps}))\le \eta\}\Big] + \E\Big[{\rm e}^{\beta\langle f,\mu_\beta\rangle}\1\{\mu_\beta(U_{\delta}(Q_{\eps}))> \eta\}\Big]\\
& \le {\rm e}^{\beta\eta K}+{\rm e}^{\beta K}\P(\mu_\beta(U_{\delta}(Q_{\eps}))>\eta)\\
&  \le {\rm e}^{\beta\eta K}+{\rm e}^{\beta K}{\rm e}^{-\beta K(1-\eta)}\le 2{\rm e}^{\beta\eta K}\le {\rm e}^{\beta K\frac{3}{4}\mu(A)}.
\end{aligned}
\end{eqnarray}
Hence, 
$$ 
J_\beta(\mu)\ge \langle f,\mu\rangle -\frac{1}{\beta}\log\E\Big[ {\rm e}^{\int_0^\beta f(B_s)\,\d s}\Big]\ge K\mu(A)-K\frac{3}{4}\mu(A)=M. 
$$
\end{proofsect}
\qed

In general, the supremum in the definition \eqref{JTKdef}  of $J_\beta$ is not attained. It is of interest to replace the function class $\Ccal_{\rm b}(\R^d)$ in \eqref{JTKdef} by some class of better behaved functions. In particular, one would like to use only functions $f$ that are extremely negative far out. This is of course possible only if $\phi^2$ decays sufficiently fast at infinity. We write $ \mu_{\beta} $ for the normalised occupation measure of the Brownian motion $ (B_s^{\ssup{1}})_{s\ge 0} $ in the following.
By $\mathfrak{m}$ we denote the initial distribution of the Brownian motions. 
\begin{lemma}[Alternate expression for $J_\beta$]\label{Jbetaident} Let $W\colon \R^d\to[0,\infty]$ be continuous in $\{W<\infty\}$, which is supposed  to contain $ \supp(\mathfrak m) $ and to be either equal to $\R^d$ or compact. Fix $\phi\in L^2(\R^d)$ satisfying $\langle W,\phi^2\rangle <\infty$.  Then
\begin{equation}\label{JWbetacompare}
J_{\beta}(\phi^2)=\sup_{h\in \Ccal_{\rm b}(\R^d)}\Big(\langle -W+h,\phi^2\rangle-\frac 1\beta\log\E\big[{\rm e}^{\beta\langle-W+h,\mu_\beta\rangle}\big]\Big).
\end{equation}
\end{lemma}

\begin{proofsect}{Proof} Let $J_{W,\beta}(\phi^2)$ denote the right-hand side of \eqref{JWbetacompare}. We first prove \lq$\geq$\rq\ in \eqref{JWbetacompare}. Let $h\in \Ccal_{\rm b}(\R^d)$. We may assume that $h\leq 0$ (otherwise we add a suitable constant to $h$). Then $f_R=(-W+h)\vee(-R)$ is a bounded continuous function for any $R>0$ with $\langle -W+h,\phi^2\rangle\leq \langle f_R,\phi^2\rangle$. Furthermore, $f_R\downarrow -W+h$ as $R\to\infty$, hence the monotonous convergence theorem yields that
$$
\lim_{R\to\infty}\E\big[{\rm e}^{\beta\langle f_R,\mu_\beta\rangle}\big]=\E\big[{\rm e}^{\beta\langle-W+h,\mu_\beta\rangle}\big].
$$
This shows that $J_{\beta}(\phi^2)\geq J_{W,\beta}(\phi^2)$ holds. Note that we did not need here that $\langle W,\phi^2\rangle <\infty$.

Now we prove  \lq$\leq$\rq\ in \eqref{JWbetacompare}. Let $f\in \Ccal_{\rm b}(\R^d)$ be given. For $R>0$, consider $h_R=(f+W)\wedge R$, then $h_R\in \Ccal_{\rm b}(\R^d)$ with $h_R\uparrow f+W$. Since $\langle W,\phi^2\rangle<\infty$, we have
$$
\liminf_{R\to\infty}\langle -W+h_R,\phi^2\rangle\ge\langle f,\phi^2\rangle.
$$
Furthermore, by the monotonous convergence theorem, 
$$
\lim_{R\to\infty}\E\big[{\rm e}^{\beta\langle -W+h_R,\mu_\beta\rangle}\big]=\E\big[{\rm e}^{\beta\langle f,\mu_\beta\rangle}\big].
$$
This implies that \lq$\leq$\rq\ holds in \eqref{JWbetacompare}.
\qed
\end{proofsect}

Let us draw a conclusion for compactly supported functions $\phi$. For a measurable set $A\subset\R^d$, we denote by $\Ccal_{\rm b}(A)$ the set of continuous bounded functions $A\to\R$.

\begin{cor} Fix $\phi\in L^2(\R^d)$ satisfying $ ||\phi||_2=1$. If the support of $\phi$ is compact, connected and contains $ \supp(\mathfrak m) $, then 
\begin{equation}
J_\beta(\phi^2)=\sup_{f\in\Ccal_{\rm b}(\supp(\phi))}\Bigl(\langle f,\phi^2\rangle-\frac{1}{\beta}\log \E\big[{\rm e}^{\beta \langle f,\mu_\beta\rangle}\1_{\{\supp (\mu_\beta)\subset\supp(\phi) \}}\big]\Bigr).
\end{equation}
\end{cor}
\begin{proofsect}{Proof} 
We pick $W=\infty\1_{\supp(\phi)}$ in Lemma~\ref{Jbetaident} and see that, on the right hand side of \eqref{JWbetacompare}, we may insert the indicator on $\{\supp (\mu_\beta)\subset\supp(\phi) \}$ in the expectation and can drop $W$ in the exponent. Hence, both this expectation and the first term, 
$\langle -W+f,\phi^2\rangle$, do  not depend on the values of $f$ outside $\supp(\phi)$.
\qed
\end{proofsect}

The next lemma shows the interplay between the arguments for the functional $ J_{\beta} $ and the fixed initial distribution of the Brownian motions.

\begin{lemma}
Let $\phi\in L^2(\R^d)$ satisfying $\|\phi\|_2=1$. If $ \dist(\supp(\phi),\supp(\mathfrak{m})) >0 $, then $ J_{\beta}(\phi^2)=+\infty $.
\end{lemma}
\begin{proofsect}{Proof}
Let $ S$ be an open neighbourhood of $ \supp (\phi) $ with $ \delta=\mbox{dist }(S,\supp(\mathfrak{m}))>0 $. Pick $K>0$ and  a continuous bounded function $ f\colon\R^d\to[0, K] $ with $ \supp (f)\subset S $ and $ f|_{\supp(\phi)}=K $. Then $ \langle f,\phi^2\rangle=K $. Then we have
\begin{equation}
\begin{aligned}
\E\big[{\rm e}^{\beta\langle f,\mu_\beta\rangle}\big]&=\E\big[{\rm e}^{\beta\langle f,\mu_\beta\rangle}\1\{\mu_\beta(S)< 1-K^{-1/2}\}\big] +\E\big[{\rm e}^{\beta\langle f,\mu_\beta\rangle}\1\{\mu_\beta(S)> 1-K^{-1/2}\}\big] \\ 
&\le {\rm e}^{\beta K(1-K^{-1/2})}+{\rm e}^{\beta K}\P\Big(\sup_{0\le t\le {K}^{-1/2}} |B_t| \ge \delta\Big)\\
&\le {\rm e}^{\beta K}\Big({\rm e}^{-\beta K^{1/2}}+\frac 2{\sqrt{2\pi}}{\rm e}^{-\frac 12\delta^2 K^{1/2}}\Big).
\end{aligned}
\end{equation}
Hence, 
$$
J_{\beta}(\phi^2)\geq \langle f,\phi^2\rangle-\frac 1\beta\log \E\big[{\rm e}^{\beta\langle f,\mu_\beta\rangle}\big]
\geq -\frac 1\beta\log\Big({\rm e}^{-\beta K^{1/2}}+\frac 2{\sqrt{2\pi}}{\rm e}^{-\frac 12\delta^2 K^{1/2}}\Big).
$$
Letting $K\to\infty$ shows that $  J_{\beta}(\phi^2)=+\infty $.
\end{proofsect}
\qed

\subsection{Analysis of $\boldsymbol {\chi^{\ssup{\otimes}}_{\sca}(\beta)}$.}\label{sec-chi}

\begin{proofsect}{Proof of Lemma~\ref{GPTKmin}}
The uniqueness of the minimiser follows from the convexity of the functionals $J_\beta $ and $\langle W,\cdot\rangle$, together with the strict convexity of $\phi^2\mapsto\|\phi\|_4^4$.

Let $(\phi_n)_{n\in\N}$ be an approximative sequence of minimisers for the formula in \eqref{GPTK}, i.e., $ \phi_n\in L^2(\R^d)$,  $\|\phi_n\|_2=1 $ for any $ n\in \N $ and 
$$
\lim_{n\to\infty}\Big(J_\beta(\phi_n^2)+\langle W,\phi_n^2\rangle+4\pi\sca\,\|\phi_n\|_4^4\Big)=\chi^{\ssup{\otimes}}_{\sca}(\beta).
$$
In particular, the sequences $(J_\beta(\phi_n^2))_n$, $(\langle W,\phi_n^2\rangle)_n$ and $(\|\phi_n\|_4)_n$ are bounded. Since $W$ explodes at infinity by Assumption (W), the sequence of probability measures $(\phi_n^2(x)\,\d x)_{n\in\N}$ is tight. According to Prohorov's Theorem, there is a probability measure $\mu$ on $\R^d$ such that $\phi_n^2(x)\,\d x$ converges weakly towards $\mu$ as $n\to\infty$, along a suitable subsequence. Since the sequence $(\phi_n)^2_n$ is bounded in $L^2(\R^d)$, the Banach-Alaoglu Theorem implies that we may assume that, along the same sequence, $(\phi_n^2)_n$ converges weakly in $L^2(\R^d)$ towards some $\phi^2\in L^2(\R^d)$. 
Since $J_\beta $ is weakly lower semi continuous (in the sense of probability measures), we get
$$
\liminf\limits_{n\to\infty}J_{\beta}(\phi_n^2)\ge J_{\beta}(\mu),
$$ and with Lemma \ref{densitiesofmeasuresJ} we conclude that $ \mu(\d x)=\widetilde \phi(x)^2 \,\d x $ for some function $ \widetilde \phi^2\in L^2(\R^d) $ with $ ||\widetilde \phi||_2=1 $. By the weak convergence in $ L^2(\R^d) $, combined with the weak convergence in the sense of probability measures, for any continuous bounded function $ \psi $ with compact support we have $ \langle \psi,\phi^2\rangle=\langle \psi,\widetilde \phi^2\rangle $. Hence, we get $ \phi^2=\widetilde \phi^2 $ a.e..
As $J_\beta $ is weakly lower semi continuous (in the sense of probability measures), $\|\cdot\|_4$ is $L^4$-weakly lower semi continuous, and $\phi\mapsto\langle W,\phi^2\rangle$ is lower semi continuous, we have that 
$$
\begin{aligned}
J_\beta(\phi^2)+\langle W,\phi^2\rangle+4\pi\sca(v)\,\|\phi\|_4^4
&\leq \liminf_{n\to\infty}\Big(J_\beta(\phi_n^2)+\langle W,\phi_n^2\rangle+4\pi\sca(v)\,\|\phi_n\|_4^4\Big)\\
&=\chi^{\ssup{\otimes}}_{(\sca)}(\beta).
\end{aligned}
$$
This shows that the limiting point $\phi$ is the minimiser in \eqref{GPTK}. This ends the proof of Lemma~\ref{GPTKmin}(i). Furthermore, the proof also shows that the minimising sequence converges, along some subsequence, towards the unique minimiser in all the three weak senses: in $L^2$, $L^4$ and weakly as a probability measures. This implies Lemma~\ref{GPTKmin}(ii).
\end{proofsect}
\qed

\section{Large-$N$ behaviour: Proof of Theorem~\ref{Kmodel,N}}\label{sec-ProofN}

In this section we prove Theorem~\ref{Kmodel,N}. We shall proceed according to the well-known G\"artner-Ellis theorem, which relates logarithmic asymptotic of probabilities to the ones of expectations of exponential integrals. Therefore, we have to establish the existence of the logarithmic moment generating function of $\overline{\mu}_{N,\beta}$ under the measure with density ${\rm e}^{-H_{N,\beta}-K^{\ssup N}_{N,\beta}}$.
The main step in the proof of Theorem~\ref{Kmodel,N} is the following.

\begin{prop}[Asymptotic for the cumulant generating function]\label{uppKN} For any $ f\in\Ccal_{\rm b}(\R^d)$ the cumulant generating function exists, i.e.,
\begin{equation}\label{logmomentfunctionKupper}
\lim_{N\to\infty}\frac 1{N\beta}\log\E\Bigl[{\rm e}^{-H_{N,\beta}-K^{\ssup N}_{N,\beta}}{\rm e}^{N\langle f,\overline{\mu}_{N,\beta}\rangle}\Bigr]= -\chi^{\ssup{\otimes}}(f),
\end{equation}
where 
\begin{equation}\label{chifdef}
\chi^{\ssup{\otimes}}(f)=\frac 1\beta\inf_{\phi\in L^2(\R^d)\colon ||\phi||_2^2=1}\Bigl\{J_\beta(\phi^2) + \langle W-f,\phi^2\rangle+4\pi{\sca}(v)||\phi||_4^4\Bigr\}.
\end{equation}
\end{prop}

Indeed, Theorem~\ref{Kmodel,N} follows from Proposition~\ref{uppKN} as follows.

\begin{proofsect}{Proof of Theorem~\ref{Kmodel,N}} We are going to use the G\"artner-Ellis Theorem, see \cite[Cor.~4.6.14]{DZ98}. For this, we only have to show that the sequence of the  $\overline{\mu}_{N,\beta}$ is exponentially tight under measure with density ${\rm e}^{-H_{N,\beta}-K^{\ssup N}_{N,\beta}}$ and that the map $f\mapsto \chi^{\ssup{\otimes}}(f)$ is G\^ateaux differentiable.

The proof of exponential tightness is easily done using our assumption that $\lim_{R\to\infty}\inf_{|x|\geq R}W(x)=\infty$ in combination with the theorems by Prohorov and Portmanteau; we omit the details.

Fix $ f\in \Ccal_{\rm b}(\R^d) $. The proof of Lemma \ref{GPTKmin} shows that the infimum in the formula of the right hand side of  \eqref{chifdef} is attained. Let $ \phi_f^2 \in L^1(\R^d) $ be the minimiser for the right hand side of \eqref{chifdef}, and, for some $ g\in\Ccal_{\rm b}(\R^d) $,  let $ \phi_{f+tg}^2\in L^1(\R^d), t > 0, $ be corresponding minimiser for $ f+tg $ instead of $ f $. We obtain
\begin{equation}\label{gateauxdiff}
\frac{1}{t}\Bigl[\chi^{\ssup{\otimes}}(f+tg)-\chi^{\ssup{\otimes}}(f)\Bigr]\ge \langle g,\phi^2_{f+tg}\rangle.
\end{equation}
Since $ (\phi^2_{f+tg})_{t\ge 0} $ is easily seen to be convergent weakly (in the sense of probability measures) towards the minimiser $ \phi^2_f $ of the right hand side of  \eqref{chifdef}, it is clear that the right hand side of \eqref{gateauxdiff} converges towards $ \langle g,\phi_f^2\rangle $. Analogously, one shows the complementary bound. This implies the G\^ateaux-differentiability of $ \chi^{\ssup{\otimes}} $ with
$$
\frac{\partial}{\partial g}\chi^{\ssup{\otimes}}(f)=\langle g,\phi_f^2\rangle.
$$ 
\qed\end{proofsect}

In Section~\ref{sec-BISLT} we give a heuristic explanation of \eqref{logmomentfunctionKupper} and introduce the Brownian intersection local times, an important object in our proof. In Sections~\ref{upperboundKNlimit} and \ref{loweboundKNlimit}, respectively, we prove the upper and the lower bound in \eqref{logmomentfunctionKupper}.

\subsection{Heuristics and Brownian intersection local times}\label{sec-BISLT}

In this section, we give a heuristic explanation of the assertion of Proposition~\ref{logmomentfunctionKupper}. We rewrite the two Hamiltonians in terms of functionals of the mean $\overline \mu_{N,\beta}$ defined in \eqref{meanmu} and use a well-known large deviation principle for $\overline\mu_{N,\beta}$. In particular, we introduce an object that will play an important role in the proofs, the Brownian intersection local times. For the definition and the most important facts on  large deviation theory used, see the Appendix or consult \cite{DZ98}.

Rewriting the first Hamiltonian in terms of $\overline \mu_{N,\beta}$ is an easy task and can be done for any fixed $N$: 
\begin{equation}\label{HVrewrite} 
H_{N,\beta}=N\beta\int_{\R^d} W(x)\frac 1N \sum_{i=1}^N \mu_\beta^{\ssup{i}}(\d x)=N\beta \bigl\langle W, \overline \mu_{N,\beta}\bigr\rangle. 
\end{equation} 
 
Now we rewrite the second Hamiltonian, which will need Brownian intersection local times and an approximation for large $N$. Let us first introduce the intersection local times, see \cite{GHR84}. For the following, we have to restrict to the case $d\in\{2,3\}$. 
 
Fix $1\leq i<j\leq N$ and consider the process $B^{\ssup{i}}-B^{\ssup{j}}$, the so-called {\em confluent Brownian motion\/} of $B^{\ssup{i}}$ and $-B^{\ssup{j}}$. This two-parameter process possesses a local time process, i.e., there is a random process $(L^{\ssup{i,j}}_\beta(x))_{x\in\R^d}$ such that, for any bounded and measurable function $f\colon\R^d\to \R$, 
\begin{eqnarray}\label{intersection}
\int_{\R^d}f(x)L^{\ssup{i,j}}_\beta(x)\,\d x=\frac 1{\beta^2}\int_0^\beta\d s\int_0^\beta\d t\, f\bigl(B_s^{\ssup{i}}-B_t^{\ssup{j}}\bigr)=\int_{\R^d}\int_{\R^d}\mu_\beta^{\ssup{i,j}}(\d x)\mu_\beta^{\ssup{i,j}}(\d y) f(x-y). 
\end{eqnarray}
Hence, we may rewrite $K^{\ssup N}_{N,\beta}$ as follows: 
\begin{equation}\label{heuristic} 
\begin{aligned} 
K^{\ssup N}_{N,\beta}&=\beta N^{d-1}\sum_{1\leq i<j\leq N} \int_{\R^d}v(z N)L_\beta^{\ssup{i,j}}(z)\,\d z\\ 
&=N\beta\int_{\R^d}v(x)\frac 1{N^2}\sum_{1\leq i<j\leq N}L_\beta^{\ssup{i,j}}({\textstyle{\frac 1N}x})\,\d x. 
\end{aligned} 
\end{equation} 
It is known \cite[Th.~1]{GHR84} that $(L_\beta^{\ssup{i,j}}(x))_{x\in \R^d} $ may be chosen continuously in the space variable. Furthermore, the random variable $L_\beta^{\ssup{i,j}}(0)=\lim_{x\to 0}L_\beta^{\ssup{i,j}}(x)$ is equal to the normalised total intersection local time of the two motions $B^{\ssup{i}}$ and $B^{\ssup{j}}$ up to time $\beta$. Formally, 
\begin{equation}\label{ISLT} 
L_\beta^{\ssup{i,j}}(0)=\frac 1{\beta^2}\int_A \d x\,\int_0^\beta \d s\,\frac{\1\{B_s^{\ssup{i}}\in \d x\}}{\d x}\int_0^\beta \d t\,\frac{\1\{B_t^{\ssup{j}}\in\d  x\}}{\d x}=\int_A \d x\,\frac{\mu_\beta^{\ssup{i}}(\d x)}{\d x}\frac{\mu_\beta^{\ssup{j}}(\d x)}{\d x}, 
\end{equation} 
Using the continuity of $L_\beta^{\ssup{i,j}}$, we approximate
$$ 
\begin{aligned}
K_{N,\beta}^{\ssup N}&\approx N\beta 4\pi\sca(v)\,\frac 2{N^2}\sum_{1\leq i<j\leq N}L_\beta^{\ssup{i,j}}(0)
\approx N\beta4\pi\sca(v)\, \Bigl\langle \frac 1N \sum_{i=1}^N \mu_\beta^{\ssup{i}},\frac 1N \sum_{i=1}^N \mu_\beta^{\ssup{i}}\Bigr\rangle\\
&= N\beta4\pi\sca(v)\,\Bigl\|\frac{\d \overline\mu_{N,\beta}}{\d x}\Bigr\|_2^2. 
\end{aligned}
$$ 
where we conceive $\mu_\beta^{\ssup{i}}$ as densities, like in \eqref{ISLT}. 

The main ingredient is now that $(\overline\mu_{N,\beta})_{N\in\N}$ satisfies a large deviation principle on $\Mcal_1(\R^d)$ with speed $N\beta$ and rate function $J_\beta$. This fact directly follows from Cram\'er's Theorem, together with the exponential tightness of the sequence $(\overline\mu_{N,\beta})_{N\in\N}$. Hence, using Varadhan's Lemma and ignoring the missing continuity of the map $\mu\mapsto \|\frac{\d \mu}{\d x}\|_2^2$, this heuristic explanation is finished by
$$
\begin{aligned}
\E\Bigl[{\rm e}^{-H_{N,\beta}-K^{\ssup N}_{N,\beta}}{\rm e}^{N\langle f,\overline{\mu}_{N,\beta}\rangle}\Bigr]
&\approx \E\Bigl[\exp\Big\{-N\beta\Big[\bigl\langle W-f, \overline \mu_{N,\beta}\bigr\rangle-4\pi\sca(v)\,\Bigl\|\frac{\d \overline\mu_{N,\beta}}{\d x}\Bigr\|_2^2\Big]\Big\}\Big]\\
&\approx{\rm e}^{-N\beta \chi^{\ssup{\otimes}}(f)},
\end{aligned}
$$
Here we substituted $\phi^2(x)\,\d x=\mu(\d x)$ and noticed that, according to Lemma~\ref{densitiesofmeasuresJ}, we may restrict the infimum over probability measures to the set of their Lebesgue densities $\phi^2$.

\subsection{Proof of the upper bound  in {\bf (\ref{logmomentfunctionKupper})} in Proposition \ref{uppKN}}\label{upperboundKNlimit}

In this section we prove the upper bound in \ref{logmomentfunctionKupper} in Proposition~\ref{uppKN}. Our proof goes along the lines of the argument sketched in Section~\ref{sec-BISLT}. However, in order to arrive at a setting in which we may apply Cram\'er's Theorem and Varadhan's Lemma, we will have to prepare with a number of technical steps. More precisely, we will have to estimate the interaction term $K_{N,\beta}^{\ssup N}$ from below in terms of a smoothed version of the intersection local times. This version will turn out to be a bounded and continuous functional of the mean of the normalised occupation times measures, $\overline\mu_{N,\beta}$, which are the central object of the analysis.

Our strategy is as follows. First, we distinguish those events on which, for at least $((1-\eta)N)^2$ pairs $(i,j)$ of indices, the intersection local times $L^{\ssup{i,j}}(x)$ for $|x|\leq 2\eps $ are sufficiently close to $L^{\ssup{i,j}}(0)$, and its complement. More precisely, we will have $|L^{\ssup{i,j}}(x)-L^{\ssup{i,j}}(0)|< \xi $ for these $(i,j)$ and $x$. Here $\xi,\eta, \eps$ are positive parameters which will eventually be sent to zero. (The complement of the event considered will turn out to be small by the continuity of the intersection local times in zero.) The replacement of $L^{\ssup{i,j}}(x)$ by $L^{\ssup{i,j}}(0)$ will require a space cutting argument, i.e., we will restrict the interaction from $\R^d$ to the cube $Q_R=[-R,R]^d$ for some $R>0$ which will eventually be sent to infinity. Our second main step is to replace $L^{\ssup{i,j}}(0)$ by the smoothed version $L^{\ssup{i,j}}\ast\kappa_\eps\ast\kappa_{\eps}(0)$, where $\kappa_{\eps}$ is a smooth approximation of the Dirac function $\delta_0$ as $\eps\downarrow 0$. The smoothed intersection local times can easily be written as a continuous bounded functional of the mean of the normalised occupation times measures, $\overline\mu_{N,\beta}$. Hence, Cram\'er's Theorem and Varadhan's Lemma become applicable, and we arrive at an upper bound for the large-$N$ rate in terms of an explicit variational formula, which depends on the parameters. Finally, we send the parameters to zero and infinity, respectively.

Let us come to the details. Introduce the following random set of index pairs,
\begin{equation}\label{Kset}
D_N=D_N(\eps,\xi)=\Big\{(i,j)\in\{1,\ldots,N\}^2\setminus\Delta_N\colon\sup_{|x|\le 2\eps}|L_\beta^{\ssup{i,j}}(0)-L_\beta^{\ssup{i,j}}(x)| <\xi\Big\},
\end{equation} 
where $ \Delta_N=\{(i,i)\colon i\in\{1,\ldots,N\}\} $ denotes the diagonal in $\{1,\ldots,N\}^2$. 
Fix $\eta>0 $ and consider the event
\begin{equation}\label{eventA}
A_N=A_N(\eps,\xi,\eta)=\big\{\exists I\subset\{1,\ldots,N\}\colon I\times I\setminus\Delta_N \subset D_N, |I|\ge (1-\eta)N\big\}
\end{equation}
In words, on $A_N$, there is a quite large set $I$ of indices such that all pairs $(i,j)$ of distinct indices in $I$ satisfy $|L_\beta^{\ssup{i,j}}(0)-L_\beta^{\ssup{i,j}}(x)| <\xi$ for all $|x|\leq 2\eps$.

First we show that the contribution coming from the complement $ A_N^{{\rm c}} $ vanishes for small $\eps$:

\begin{lemma}[$ A_N^{{\rm c}}$ is negligible]\label{eventAcompl}
For any $ \xi>0$ and any $\eta \in(0,\frac 12) $,
\begin{equation}\label{AcompsmallHK}
\lim_{{\eps}\to 0}\limsup_{N\to\infty}\frac{1}{N}\log \E\Big[{\rm e}^{-H_{N,\beta}-K^{\ssup N}_{N,\beta}}{\rm e}^{N\langle f,\overline{\mu}_{N,\beta}\rangle}\1_{A_N^{{\rm c}}(\eps,\xi,\eta)}\Big]=-\infty.
\end{equation}
\end{lemma}

\begin{proofsect}{Proof}
Since $H_{N,\beta}$, $K_{N,\beta}^{\ssup N}$ and $f$  are bounded from below, it suffices to show that
\begin{equation}\label{Acompsmall}
\lim_{\eps\to 0}\limsup_{N\to\infty}\frac{1}{N}\log \P\big(A_N^{{\rm c}}(\eps,\xi,\eta)\big)=-\infty.
\end{equation}
Note that
\begin{equation}
A_N^{{\rm c}}=A_N^{{\rm c}}(\eps,\xi,\eta)=\{\forall I\subset\{1,\ldots,N\}\colon |I|\ge (1-
\eta)N\Longrightarrow (I\times I)\setminus\Delta_N \not\subset D_N\}.
\end{equation}
First we show that, on $ A_N^{{\rm c}} $, there are pairwise different integers $ i_1,j_1,i_2,j_2,\ldots,i_{\lfloor\eta N/3\rfloor},j_{\lfloor\eta N/3\rfloor} $ in $ \{1,\ldots,N\} $ such that $ (i_l,j_l)\in D_N^{{\rm c}} $ for all $ l=1,\ldots,\lfloor\eta N/3\rfloor$.

We construct the indices inductively. Consider $ I_1=\{1,\ldots, \lceil(1-\eta)N\rceil\} $, then
$$
I_1\times I_1\setminus \Delta_N\not\subset D_N,
$$ 
i.e., there is a pair $ (i_1,j_1)\in (I_1\times I_1\setminus\Delta_N)\cap D_N^{{\rm c}}$. Also the set
$$
I_2=(I_1\setminus\{i_1,j_1\})\cup\{(1- \eta)N+1,(1- \eta)N+2\}
$$ 
has no less than $ (1- \eta)N $ elements. Hence there is a pair $ (i_2,j_2)\in (I_2\times I_2\setminus\Delta_N)\cap D_N^{{\rm c}} $. Clearly $ \sharp\{i_1,j_1,i_2,i_2\}=4 $. In this way, we can proceed altogether at least $ \floor{\frac{1}{2}\eta N} $ times. This procedure constructs the indices $ i_1,j_1,i_2,j_2,\ldots,i_{ \lfloor\eta N/3\rfloor},j_{\lfloor\eta N/3\rfloor} $ pairwise different such that $ (i_1,j_1),\ldots,(j_{\lfloor\eta N/3\rfloor},j_{\lfloor\eta N/3\rfloor})\in D_N^{{\rm c}} $.

Now we prove that \eqref{Acompsmall} holds. We abbreviate \lq p.d.\rq\ for `pairwise disjoint` in the following. For notational convenience, we drop the brackets $\lfloor\cdot\rfloor$. Because of the preceding, we have
\begin{equation}\label{uppercompl}
\begin{aligned}
\P(A_N^{{\rm c}}(\eps,\xi,\eta))
&\le \sum_{i_1,j_1,\ldots,i_{\eta N/3},j_{\eta N/3}\in \{1,\ldots, N\}, \mbox{{\tiny  p.d.}}}\P\Bigl(\forall l=1,\ldots,\eta N/3\colon\sup_{|x|\le 2\eps}|L_\beta^{\ssup{i_l,j_l}}(x)-L_\beta^{\ssup{i_l,j_l}}(0)|\ge \xi\Bigr)\\
&=\left(\begin{array}{c}N\\2\eta N/3\end{array}\right)\P\Bigl(\sup_{|x|\le 2\eps}|L_\beta^{\ssup{1,2}}(x)-L_\beta^{\ssup{1,2}}(0)|\ge \xi\Bigr)^{\eta N/3}\\ &\le \exp\Bigl(-N\Bigl[\log\frac{1}{2}-\frac\eta3\log\P\Bigl(\sup_{|x|\le 2\eps}|L_\beta^{\ssup{1,2}}(x)-L_\beta^{\ssup{1,2}}(0)|\ge \xi\Bigr)\Bigr]\Bigr).
\end{aligned}
\end{equation}  
Since the process $(L_\beta^{\ssup{1,2}}(x))_{x\in\R^d}$ may be chosen continuously in the space variable \cite[Th.~1]{GHR84}, we have
\begin{equation}
\lim_{\eps\to 0} \P\Bigl(\sup_{|x|\le 2\eps}|L_\beta^{\ssup{1,2}}(x)-L_\beta^{\ssup{1,2}}(0)|\ge \xi\Bigr)=0.
\end{equation} 
This, together with (\ref{uppercompl}), concludes the proof.
\end{proofsect}
\qed

Now we estimate $K^{\ssup N}_{N,\beta}$ on the event $A_N$. For $R>0$, we recall that $Q_R=[-R,R]^d$ and introduce 
\begin{equation}
{\sca}_R(v)=\frac 1{8\pi}\int_{Q_R}v(|x|)\,\d x.
\end{equation}
Let $\kappa\colon\R^d\to[0,\infty)$ be a smooth function with support in $[-1,1]^d$ and $\int\kappa(x)\,\d x=1$. For $\eps>0$, we define $\kappa_{\eps}(x)=\eps^{-d}\kappa(x/\eps)$. Then $\kappa_\eps$ is an approximation of $\delta_0$ as $\eps\downarrow 0$, and we have $ \mbox{supp }\kappa_{\eps}\subset Q_{\eps}$ and $ \int_{\R^d}\kappa_{\eps}(x)\,\d x=1 $ for any $\eps>0$.

\begin{lemma}[Estimating $K^{\ssup N}_{N,\beta}$ on $A_N(\eps,\xi,\eta)$]\label{lem-Kesti} Fix $\eps,\xi,\eta>0$. Then, for any  $R>0$ and any $N\in\N$ satisfying $N>R/(2\eps)$, on the event $A_N(\eps,\xi,\eta)$,
\begin{equation}\label{estimationUpper4}
-K^{\ssup N}_{N,\beta}  \le -4\pi{\sca}_R(v)\beta |I|\,(1- \eta)\big\|\overline{\mu}_{I,\beta}\ast\kappa_{\eps}\big\|_2^2 
+8\pi\xi{\sca}_R(v)+8\pi\frac{{\sca}_R(v)||\kappa_{\eps}||_{\infty}}{N},
\end{equation} 
where the random subset $I$ of $\{1,\dots,N\}$ in \eqref{eventA} is chosen minimally with $|I|\geq (1- \eta)N$ and $(I\times I\setminus \Delta_N)\subset D_N$, and
\begin{equation}
\overline{\mu}_{I,\beta}=\frac{1}{|I|}\sum_{i\in I}\mu_\beta^{\ssup{i}}
\end{equation} 
denotes the mean of the corresponding normalised occupation measures.
\end{lemma}

\begin{proofsect}{Proof}
First, we write the interaction terms for the scaled pair potential $ v_N $ as integrals against the intersection local times of two Brownian motions at spatial points $ x/N $. As the pair interaction $ v $ is positive we get easily an upper bound, when we restrict the integrations to the box $ Q_R $. On this compact set we use the continuity to get rid of the dependence of the spatial variables in the integrals. We restrict the summation over all pairs of indices to the random set $D_N=D_N(\eps,\xi) $. 

Recall \eqref{heuristic}. On the event  $ A_N(\eps,\xi,\eta) $, we may estimate, for all $ N>R/(2\eps) $,
\begin{eqnarray}\label{estimationUpper1}
\begin{aligned}
-\frac 1{N\beta}K^{\ssup N}_{N,\beta}&=-\frac{1}{N^2}\sum_{1\le i <j\le N}\int_{\R^d}v(|x|)L_\beta^{\ssup{i,j}}\Big(\frac{x}{N}\Big)\,\d x
\le -\frac{1}{N^2}\sum_{(i,j)\in D_N}\int_{Q_R}v(|x|)L_\beta^{\ssup{i,j}}\Big(\frac{x}{N}\Big)\,\d x \\
& \le -4\pi\frac{{\sca}_R(v)}{N^2}\sum_{(i,j)\in D_N} L_\beta^{\ssup{i,j}}(0) + \frac{1}{N^2}\sum_{(i,j)\in D_N}\int_{Q_R} v(|x|)\Bigl|L_\beta^{\ssup{i,j}}(0) -L_\beta^{\ssup{i,j}}(x/N)\Bigr|\d x\\ 
& \le -4\pi\frac{{\sca}_R(v)}{N^2}\sum_{(i,j)\in D_N} L_\beta^{\ssup{i,j}}(0) +\frac{\xi |D_N|^2}{N^2}8\pi{\sca}_R(v)\\
&\le -4\pi\frac{{\sca}_R(v)}{N^2}\sum_{(i,j)\in D_N} L_\beta^{\ssup{i,j}}(0) +8\pi{\xi}{\sca}_R(v).
\end{aligned}
\end{eqnarray} 

For any $ (i,j)\in D_N $, we now replace the intersection local time at zero, $L_\beta^{\ssup{i,j}}(0)$, with the smoothed version $ L_\beta^{\ssup{i,j}}\ast\kappa_{\eps}\ast\kappa_{\eps}(0) $. The replacement error is estimated by
\begin{equation}\label{errorsmothing}
\Bigl|L_\beta^{\ssup{i,j}}(0)-L_\beta^{\ssup{i,j}}\ast\kappa_{\eps}\ast\kappa_{\eps}(0)\Bigr|\le \int_{Q_{2\eps}}\d x\, \kappa_{\eps}\ast\kappa_{\eps}(x)\Bigl|L_\beta^{\ssup{i,j}}(0)-L_\beta^{\ssup{i,j}}(x)\Bigr|\le \xi,\qquad (i,j)\in D_N.
\end{equation}
Hence, we may continue the estimation in (\ref{estimationUpper1}) by
\begin{equation}\label{estimationUpper2}
-\frac 1{N\beta}K^{\ssup N}_{N,\beta} \le -4\pi\frac{{\sca}_R(v)}{N^2}\sum_{(i,j)\in D_N} L_\beta^{\ssup{i,j}}\ast\kappa_{\eps}\ast\kappa_{\eps}(0)  + 8\pi\xi{\sca}_R(v).
\end{equation}
Recall the defining property on the intersection local time in \eqref{intersection}. Hence, we can write the smoothed version in terms of a convolution of the normalised occupation measures as follows:
\begin{equation}\label{intersec-occ}
\begin{aligned}
L_\beta^{\ssup{i,j}}\ast\kappa_{\eps}\ast\kappa_{\eps} (0)&=\int_{\R^d}\int_{\R^d}\int_{\R^d}\mu_\beta^{\ssup{i}}(\d x)\mu_\beta^{\ssup{j}}(\d y)\d z\,\kappa_{\eps}(y-x-z)\kappa_{\eps}(z)\\&=\int_{\R^d}\int_{\R^d}\int_{\R^d}\mu_\beta^{\ssup{i}}(\d x)\mu_\beta^{\ssup{j}}(\d y)\d w\,\kappa_{\eps}(w-y)\kappa_{\eps}(w-x)\\
&=\big\langle \mu_\beta^{\ssup i}\ast\kappa_\eps,\mu_\beta^{\ssup j}\ast\kappa_\eps\big\rangle.
\end{aligned}
\end{equation}

To proceed with our estimates from (\ref{estimationUpper2}), we add now the self-intersection terms, i.e., the diagonal terms where $i=j$. The additional terms are bounded by the $ L^{\infty}$-norm of $ \kappa_{\eps} $ as seen from
\begin{equation}
\begin{aligned}
\big\langle \mu_\beta^{\ssup i}\ast\kappa_\eps,\mu_\beta^{\ssup j}\ast\kappa_\eps\big\rangle
&=\int_{\R^d}\int_{\R^d}\int_{\R^d}\d x\,\mu_\beta^{\ssup{i}}(\d y)\,\mu_\beta^{\ssup{j}}(\d z)\,\kappa_{\eps}(x-y)\kappa_{\eps}(x-z)\\
&\le ||\kappa_{\eps}||_{\infty}\int_{\R^d}\int_{\R^d}\d x\,\mu_\beta^{\ssup{i}}(\d y)\,\kappa_{\eps}(x-y)=||\kappa_{\eps}||_{\infty}.
\end{aligned}
\end{equation} 
Hence, we obtain from \eqref{estimationUpper2} that
\begin{equation}\label{estimationUpper3}
\begin{aligned}
-\frac 1{N\beta}K^{\ssup N}_{N,\beta}& \le -4\pi\frac{{\sca}_R(v)}{N^2}\sum_{(i,j)\in D_N\cup \Delta_N}\big\langle\mu_\beta^{\ssup{i}}\ast\kappa_{\eps},\m_\beta^{\ssup{j}}\ast\kappa_{\eps}\big\rangle +8\pi\xi{\sca}_R(v)+4\pi\frac{{\sca}_R(v)||\kappa_{\eps}||_{\infty}}{N}.
\end{aligned}
\end{equation}

Recall that we work on the event $A_N$ defined in \eqref{eventA}. On this event, $D_N\cup\Delta_N$ contains a subset of the form $I\times I$ with $I\subset\{1,\dots,N\}$ and $|I|\ge (1- \eta)N $. Let such a random set be chosen, for definiteness we choose it minimally. We continue the estimation of (\ref{estimationUpper3}) with
\begin{equation}
-\frac 1{N\beta}K^{\ssup N}_{N,\beta}  \le -4\pi\frac{{\sca}_R(v)|I|^2}{N^2}\big\langle\overline{\mu}_{I,\beta}\ast\kappa_{\eps},\overline{\mu}_{I,\beta}\ast\kappa_{\eps}\big\rangle 
+8\pi\xi{\sca}_R(v)+4\pi\frac{{\sca}_R(v)||\kappa_{\eps}||_{\infty}}{N},
\end{equation} 
and from this the assertion follows.
\qed
\end{proofsect}

Using Lemma~\ref{lem-Kesti} on the left hand side of \eqref{logmomentfunctionKupper}, we obtain the following bound.

\begin{cor}Fix $\eps,\xi,\eta>0$. Then, for any $R>0$, as $N\to\infty$,
\begin{equation}\label{estionA}
\begin{aligned}
\E\Bigl[&{\rm e}^{-H_{N,\beta}-K^{\ssup N}_{N,\beta}}{\rm e}^{N\langle f,\overline{\mu}_{N,\beta}\rangle}\1_{A_N(\eps,\xi,\eta)}\Bigr]
\leq {\rm e}^{o(N)} {\rm e}^{N[\eta \beta||f||_{\infty}+8\xi \beta{\sca}_R(v)+C\eta]}\\
&\times \E\Bigl[\exp\Bigl\{-\beta\lfloor (1-\eta)N\rfloor\Bigl(\langle W-f,\overline{\mu}_{\lfloor (1-\eta)N\rfloor,\beta}\rangle+4\pi{\sca}_R(v)(1- \eta)\big\|\overline{\mu}_{\lfloor (1-\eta)N\rfloor,\beta}\ast \kappa_{\eps}\big\|_2^2\Bigr)\Bigr\}\Bigr],
\end{aligned}
\end{equation}
where $C>0$ is an absolute constant. 
\end{cor}

\begin{proofsect}{Proof}
Since the trapping potential $ W $ is nonnegative, we easily estimate
\begin{equation}
-H_{N,\beta}=-\beta\sum_{i=1}^N\langle W,\mu_\beta^{\ssup{i}}\rangle\le -\beta \sum_{i\in I}\langle W,\mu_\beta^{\ssup{i}}\rangle=-\beta|I|\langle W,\overline \mu_{I,\beta}\rangle.
\end{equation} 
Furthermore, it is easy to see that
\begin{equation}
N\beta\langle f, \overline\mu_{N,\beta}\rangle\le \beta |I|\,\langle f, \overline{\mu}_{I,\beta}\rangle + N\eta \beta\,||f||_{\infty}.
\end{equation}
Using these two estimates and \eqref{estimationUpper4} in the expectation on the left hand side of \eqref{logmomentfunctionKupper}, we obtain
\begin{equation}\label{estimationUpper5}
\begin{aligned}
\E\Bigl[&{\rm e}^{-H_{N,\beta}-K^{\ssup N}_{N,\beta}}{\rm e}^{N\langle f,\overline{\mu}_{N,\beta}\rangle}\1_{A_N(\eps,\xi,\eta)}\Bigr]\\
& \le \E\Bigl[\exp\Bigl\{-\beta|I|\Bigl(\langle W-f,\overline{\mu}_{I,\beta}\rangle+4\pi{\sca}_R(v) (1- \eta)\big\|\overline{\mu}_{I,\beta}\ast\kappa_{\eps}\big\|_2^2\Bigr)\Bigr\} \1_{A_N(\eps,\xi,\eta)}\Bigr]{\rm e}^{NC_{\eta,\xi,R}+o(N)},
\end{aligned}
\end{equation}
where
\begin{equation}
C_{\eta,\xi,R}= \eta \beta||f||_{\infty}+8\pi\xi \beta{\sca}_R(v).
\end{equation}
Now we sum over all possible values of the random set $ I $ and note that the distribution of $\overline{\mu}_{I,\beta}$ is equal to the one of $ \overline{\mu}_{l,\beta}=\frac{1}{l}\sum_{i=1}^l\mu_\beta^{\ssup{i}} $. Hence we get
\begin{equation}\label{sumestimationl}
\begin{aligned}
\mbox{l.h.s. of (\ref{estimationUpper5})}  &\le {\rm e}^{N C_{\eta,\xi,R}+o(N)} \sum_{l=\lceil (1- \eta)N\rceil}\sum_{L\subset\{1,\ldots,N\}\colon |L|=l}\\ 
&\qquad \E\Big[\1_{I=L} \exp\Big\{-l\beta\Big(\langle W-f,\overline{\mu}_{I,\beta}\rangle+4\pi{\sca}_R(v) (1- \eta)\big\|\overline{\mu}_{L,\beta}\ast\kappa_{\eps}\big\|_2^2\Big)\Big\} 
\1_{A_N(\eps,\xi,\eta)} \Big]\\
&\le {\rm e}^{NC_{\eta,\xi,R}+o(N)}\sum_{l=\lfloor (1- \eta)N\rfloor}^N\binom Nl\\ 
& \qquad\times\E\Big[\exp\Big\{-l\beta\Big(\langle W-f,\overline{\mu}_{l,\beta}\rangle+4\pi{\sca}_R(v)(1- \eta)\big\|\overline{\mu}_{l,\beta}\ast\kappa_{\eps}\big\|_2^2\Big)\Big\}\Big],
\end{aligned}
\end{equation} 
It is clear that there is a $C>0$ such that, for any $\eta>0$ and any $l\in\{\lceil (1- \eta)N\rceil,\dots,N\}$ we can estimate $\binom Nl\leq {\rm e}^{C\eta N}$. In the exponent, estimate $l\geq \lfloor(1- \eta)N\rfloor$ where $l$ is multiplied with a nonnegative factor, respectively estimate $l\langle f,\overline{\mu}_{l,\beta}\rangle \leq \lfloor(1- \eta)N\rfloor\langle f,\overline{\mu}_{l,\beta}\rangle+ \eta N\|f\|_\infty$. The sum on $l$ is estimated against $N$, which is absorbed in the term ${\rm e}^{o(N)}$.
\qed
\end{proofsect} 

Now we use arguments from large deviation theory to identify the large-$N$ rate of the right hand side of \eqref{estionA}:

\begin{lemma}[Large deviation rate]\label{Aesti} For any $\eps,\xi,\eta>0$ and any $R>0$,
\begin{equation}\label{estionAneu}
\begin{aligned}
\limsup_{N\to\infty}&\frac 1{N\beta}\log\E\Bigl[{\rm e}^{-H_{N,\beta}-K^{\ssup N}_{N,\beta}}{\rm e}^{N\langle f,\overline{\mu}_{N,\beta}\rangle}\1_{A_N(\eps,\xi,\eta)}\Bigr]\\
&\leq  \eta \beta||f||_{\infty}+8\pi\xi \beta{\sca}_R(v)+C\eta-(1- \eta)\chi^{\ssup{\otimes}}_{\eps,R,\eta}(f),
\end{aligned}
\end{equation}
where
\begin{equation}\label{chifetaRdef}
\begin{aligned}
\chi^{\ssup{\otimes}}_{\eps,R,\eta}(f)=\frac 1\beta\inf_{\phi\in L^2(\R^d)\colon \|\phi\|_2^2=1}\Bigl\{J_\beta(\phi^2)+\langle W-f,\phi^2\rangle +4\pi{\sca}_R(v)(1- \eta)\|\phi^2\ast\kappa_{\eps}\|_2^2\Bigr\}.
\end{aligned}
\end{equation}
\end{lemma}

\begin{proofsect}{Proof}
On the right hand side of \eqref{estionA}, we will apply  Cram\'er's Theorem \cite[Th.~6.1.3]{DZ98} in combination with Varadhan's Lemma \cite[Lemma~4.3.6]{DZ98} for the mean $ \overline{\mu}_{l,\beta} $ of normalised occupation measures. Let us explain what these theorems say and how we apply them. See the Appendix, Section~\ref{LDP} for a brief account on large deviation theory.

Note that $ \Mcal_1(\R^d) $ is a closed convex subset of the  space of all finite signed measures, $ \Mcal(\R^d) $. The duality relation
\begin{equation}
(f,\nu)\in \Ccal_{\rm b}(\R^d)\times\Mcal(\R^d)\mapsto \int_{\R^d}f(x)\, \nu(\d x)
\end{equation} 
determines a representation of $ \Mcal(\R^d)^{\ast} $, the topological dual of $ \Mcal(\R^d) $, as $ \Ccal_{\rm b}(\R^d) $. 
Clearly, the topology inherited by $ \Mcal_1(\R^d) $ from  $ \Mcal(\R^d) $ is the weak topology, which is induced by the test integrals against all  bounded continuous  functions. $ \Mcal_1(\R^d) $ is a Polish space  with the L\'{e}vy metric \cite{DS01}. Thus, all the assumptions of \cite[Th.~6.1.3]{DZ98} are satisfied. Hence, $(\overline \mu_{l,\beta})_{l\in\N}$ satisfies a weak large deviation principle. The rate function is equal to the Legendre-Fenchel transform of the logarithmic moment generating function of $ \mu_\beta^{\ssup{1}} $, that is, it is the function $J_\beta$ defined in \eqref{JTKdef}. 

For any $M>0$, the set 
$$
F_M=\{\mu\in\Mcal_1(\R^d)\colon \langle W,\mu\rangle\leq M\}
$$ 
is compact, as is easily derived with the help of Prohorov's Theorem, using that $\lim_{|x|\to\infty}W(x)=\infty$ and the lower semi continuity of the map $\mu\mapsto \langle W,\mu\rangle$ (see Assumption~(W)). 

We abbreviate $l=\lfloor(1-\eta)N\rfloor$ for a while. On the right hand side of \eqref{estionA}, we insert $\1_{F_M}(\overline\mu_{l,\beta})+\1_{F_M^{\rm c}}(\overline\mu_{l,\beta})$. On the event $\{\overline\mu_{l,\beta}\in F_M^{\rm c}\}$, we estimate the trap part, $\langle W,\overline\mu_{l,\beta}\rangle$, from below  against $M$ and use that the remaining terms in the exponential are bounded from below. Hence,
\begin{equation}\label{LDPappl}
\begin{aligned}
\E\Bigl[&\exp\Bigl\{-l\beta\Bigl(\langle W-f,\overline{\mu}_{l,\beta}\rangle+4\pi{\sca}_R(v)(1- \eta)\big\|\overline{\mu}_{l,\beta}\ast \kappa_{\eps}\big\|_2^2\Bigr)\Bigr\}\Bigr]\\
&\leq {\rm e}^{-l\beta(M-\|f\|_\infty)}+\E\Bigl[\exp\Bigl\{-l\beta\Bigl(\langle W-f,\overline{\mu}_{l,\beta}\rangle+4\pi{\sca}_R(v)(1- \eta)\big\|\overline{\mu}_{l,\beta}\ast \kappa_{\eps}\big\|_2^2\Bigr)\Bigr\}\1_{F_M}(\overline\mu_{l,\beta})\Bigr].
\end{aligned}
\end{equation}
The functional $\mu\mapsto \langle W-f,\mu\rangle+4\pi{\sca}_R(v)(1- \eta)\|\mu\ast \kappa_{\eps}\|_2^2$ is lower semi continuous, as is easily seen from Fatou's lemma, and bounded from below. Furthermore, according to \cite[Th.~4.5.3]{DZ98}, the family $(\overline\mu_{l,\beta})_{l\in\N}$ satisfies the upper bound in the weak large deviation principle with rate function $J_\beta$. Hence, we may apply the upper-bound part of Varadhan's lemma \cite[Lemma~4.3.6]{DZ98} to the right-hand side of \eqref{LDPappl}, to obtain, if $M$ is sufficiently large, 
\begin{eqnarray}
\begin{aligned}
\limsup_{l\to\infty}\frac{1}{l\beta}&\log\E\Bigl[\exp\Bigl\{-l\beta\Bigl(\langle W-f,\overline{\mu}_{l,\beta}\rangle+4\pi{\sca}_R(v)(1- \eta)\big\|\overline{\mu}_{l,\beta}\ast \kappa_{\eps}\big\|_2^2\Bigr)\Bigr\}\Bigr]\\
& \le -\inf_{\nu\in F_M}\Bigl(J_{\beta}(\nu)+\langle W-f,\nu\rangle +4\pi{\sca}_R(v)(1-\eta)\big\|\nu\ast\kappa_{\eps}\big\|_2^2\Bigr).
\end{aligned}
\end{eqnarray}
Lemma~\ref{densitiesofmeasuresJ} implies that the infimum on the right hand side is equal to $\chi^{\ssup{\otimes}}_{\eps,R,\eta}(f)$.
\qed
\end{proofsect}

Summarising the contributions from $A_N$ and $A_N^{\rm c}$, we obtain the following.

\begin{cor}
\begin{equation}\label{estionAneu2}
\begin{aligned}
\limsup_{N\to\infty}\frac 1{N\beta}\log\E\Bigl[{\rm e}^{-H_{N,\beta}-K^{\ssup N}_{N,\beta}}{\rm e}^{N\langle f,\overline{\mu}_{N,\beta}\rangle}\Bigr]\leq -\liminf_{\eta\downarrow 0}\liminf_{R\to \infty}\liminf_{\eps\downarrow 0}\chi^{\ssup{\otimes}}_{\eps,R,\eta}(f),
\end{aligned}
\end{equation}
where $\chi^{\ssup{\otimes}}_{\eps,R,\eta}(f)$ is defined in \eqref{chifetaRdef}.
\end{cor}

\begin{proofsect}{Proof}
A combination of Lemmas~\ref{eventAcompl} and \ref{Aesti} gives that the left hand side of \eqref{estionAneu2} is not smaller than
$$
\eta \beta||f||_{\infty}+\xi \beta{\sca}_R(v)+C\eta-(1- \eta)\liminf_{\eps\downarrow 0}\chi^{\ssup{\otimes}}_{\eps,R,\eta}(f),
$$
for any $\eta,R,\xi>0$. Letting $\eta,\xi\downarrow0$ and $R\to\infty$, we arrive at the assertion.
\qed\end{proofsect}

Now we identify the right hand side of \eqref{estionAneu2}:

\begin{lemma}[Approximating the variational formula]\label{lem-ApprVP}
\begin{equation}
\liminf_{\eta\downarrow 0}\liminf_{R\to \infty}\liminf_{\eps\downarrow 0}\chi^{\ssup{\otimes}}_{\eps,R,\eta}(f)
\geq \chi^{\ssup{\otimes}}(f),
\end{equation}
where $\chi^{\ssup{\otimes}}(f)$ is defined in \eqref{logmomentfunctionKupper}.
\end{lemma}

\begin{proofsect}{Proof} We first fix $\eta>0$ and $R>$ and prove that
\begin{eqnarray}\label{parameter4}
\liminf_{\eps\downarrow 0}\chi^{\ssup{\otimes}}_{\eps,R,\eta}(f)\ge \chi_{R,\eta}^{\ssup{\otimes}},
\end{eqnarray} 
with
\begin{equation}\label{chiRetafdef}
\chi^{\ssup{\otimes}}_{R,\eta}(f)=\frac 1\beta\inf_{\phi\in L^2(\R^d)\colon\|\phi\|_2=1}\Bigl(J_\beta(\phi^2)+\langle W-f,\phi^2\rangle +4\pi{\sca}_R(v)(1-\eta)\,||\phi||_4^4\Bigr).
\end{equation}
Let $ (\phi_{\eps})_{\eps\ge 0} $ be an approximate minimising sequence for the right hand side of \eqref{chiRetafdef}, i.e.,
$$
\liminf\limits_{\eps\searrow 0}\frac{1}{\beta}\Bigl(J_{\beta}(\phi_{\eps}^2)+\langle W-f,\phi_{\eps}^2\rangle +4\pi\sca_R(v)(1-\eta)||\phi_{\eps}^2*\kappa_{\eps}||^2_2\Bigr)=\liminf\limits_{\eps\searrow 0}\chi_{\eps,R,\eta}^{\ssup{\otimes}}(f).
$$
In particular, $ (J_{\beta}(\phi_{\eps}^2))_{\eps\ge 0}, (\langle W-f,\phi_{\eps}\rangle)_{\eps\ge 0} $ and $ (||\phi_{\eps}^2*\kappa_{\eps}||^2_2)_{\eps\ge 0} $ are bounded.
As $\eps\downarrow 0$, along suitable subsequences, the probability measures $ \mu_\eps(\d x)=\phi_\eps^2(x)\,\d x $ converge weakly to a probability measure $ \mu(\d x) $. Certainly, we may assume that $\chi^{\ssup{\otimes}}_{\eps,R,\eta}(f)$ is bounded as $\eps\downarrow 0$, and therefore also $J_\beta(\phi_{\eps}^2)$ is. The lower semi continuity of $J_\beta$ with respect to the weak topology of probability measures gives that $J_\beta(\mu)\le \liminf_{\eps\searrow 0} J_\beta(\mu_{\eps}) <\infty$.
From Lemma \ref{densitiesofmeasuresJ} we get the existence of a density for the measure $ \mu $, i.e., $ \mu(\d x)=\phi^2(x)\d x $ for some $ \phi^2\in L^2(\R^d) $ satisfying $\|\phi\|_2=1 $. Since $ \eps\mapsto ||\phi^2_{\eps}*\kappa_{\eps}||_2 $ is bounded, there is a function $ \psi\in L^4(\R^d) $ such that, along suitable subsequences, $ \phi_{\eps}^2*\kappa_{\eps} $ converges weakly in $ L^2(\R^d) $ to $ \psi^2 $, and
$ ||\psi||_4^4\le \liminf_{\eps\downarrow 0}||\phi_{\eps}^2*\kappa_{\eps}||_2^2 $. Thus, for every $ g\in \Ccal_{\rm c}(\R^d)  $ we get
\begin{eqnarray}\label{estdifferencePsiPhi}
\begin{aligned}
|\langle g,\phi^2-\psi^2\rangle| \le |\langle g,\phi^2\rangle -\langle g,\phi_{\eps}^2\rangle |+|\langle g,\phi_{\eps}^2\rangle-\langle g,\phi_{\eps}^2*\kappa_{\eps}\rangle|+|\langle g,\phi_{\eps}^2*\kappa_{\eps}-\psi^2\rangle|.  
\end{aligned}
\end{eqnarray}
Now, the first term on the right hand side of (\ref{estdifferencePsiPhi}) vanishes in the limit $ \eps\to 0 $ because $ g $ is bounded and continuous and the convergence follows from the weak convergence of the probability measures. The second term on the right hand side of (\ref{estdifferencePsiPhi}) is estimated as
\begin{equation}
|\langle g,\phi_{\eps}^2\rangle-\langle g,\phi_{\eps}^2*\kappa_{\eps}\rangle|=|\langle g-g*\kappa_{\eps},\phi_{\eps}^2\rangle| \le || g-g*\kappa_{\eps}||_{\infty}\to 0\,\mbox{ for } \eps\to 0,
\end{equation} 
since $ g\in\Ccal_{\rm c}(\R^d) $. The last term on the right hand side of (\ref{estdifferencePsiPhi}) vanishes in the limit $ \eps\to 0 $ because $ g\in L^2(\R^d) $.
Hence, $ \langle g,\phi^2\rangle=\langle g,\psi^2\rangle $ and therefore $ \phi=\psi $ almost everywhere. Clearly, $ \liminf_{\eps\downarrow 0}\langle W,\phi_{\eps}^2\rangle\ge \langle W,\phi^2\rangle. $ Altogether, (\ref{parameter4}) follows.

We now finish the proof of the lemma by showing that  
\begin{equation}\label{parameter1}
\liminf_{\eta \downarrow 0}\liminf_{R \to\infty} \chi^{\ssup{\otimes}}_{R,\eta}(f)  \ge \chi^{\ssup{\otimes}}(f).
\end{equation} 
Since the map $x\mapsto v(|x|)$ is assumed integrable, it is clear that, as $R\to\infty$,  $\sca_R(v)$ converges towards $\sca(v)$ defined in \eqref{vnorm}. Hence, the proof of \eqref{parameter1} is an easy task, and we omit it.
\qed\end{proofsect}

\subsection{Proof of the lower bound of {\bf (\ref{logmomentfunctionKupper})} in Proposition \ref{uppKN}}\label{loweboundKNlimit}

Now we turn to the proof of the lower bound in \eqref{logmomentfunctionKupper}. We were not able to produce a proof along the lines of usual large deviation arguments including Cram\'er's Theorem and Varadhan's Lemma, since we did not find any way to overcome the technical difficulties stemming from the singularity of the pair interaction term. Instead, we write the expectation on the left-hand side of \eqref{logmomentfunctionKupper} as $N$ iterated expectations over the $N$ Brownian motions and use a lower estimate that is directly implied by the definition of the rate function, $J_\beta$, more precisely, of some variant to be introduced below. In order to explain this idea, fix $\phi\in L^2(\R^d)$ and consider the random potential
\begin{equation}\label{randompot}
q_j:=-\sum_{i=1}^{j-1}V_N *\mu_\beta^{\ssup{i}} -(N-j)V_N * \phi^2,\qquad j=1,\ldots,N.
\end{equation} 
Here we used the notation $V*\mu(x)=\int_{\R^d} v(|x-y|)\,\mu(\d y)$ (analogously with $\phi^2$, conceived as a finite measure) and an analogous notation for $V_N$ with $v$ replaced by $v_N(\cdot)=N^{d-1}v(\,\cdot\, N)$.
We rewrite the left-hand side of \eqref{logmomentfunctionKupper} as follows. We write $\E^{\ssup i}$ for the expectation with respect to the $i$-th Brownian motion. Recalling the definition of $K_{N,\beta}^{\ssup N}$ in \eqref{KNbetaNdef} and noting that $q_N=-\sum_{i<N}V_N *\mu_\beta^{\ssup{i}}$, we see that we have
\begin{equation}\label{estimationstart}
\begin{aligned}
\E\Big[{\rm e}^{-H_{N,\beta}-K_{N,\beta}^{\ssup N}+\beta N\langle f,\overline{\mu}_{N,\beta}\rangle}\Big] = \E^{\ssup{1}}\otimes\cdots\otimes\E^{\ssup{N-1}}\Big[&{\rm e}^{-H_{N-1,\beta}-K_{N-1,\beta}^{\ssup N}+\beta (N-1)\langle f,\overline{\mu}_{N-1,\beta}\rangle}\\ 
&\qquad\times \E^{\ssup{N}}\big[{\rm e}^{\beta\langle q_N-W+f,\mu_\beta^{\ssup{N}}\rangle}\big]\Big].
\end{aligned}
\end{equation}
The main idea in our proof of the lower bound is that the definition of $J_\beta$ directly implies the estimate
\begin{equation}\label{toestimate}
\E\big[{\rm e}^{\beta\langle h,\mu_\beta\rangle}\big]\ge {\rm e}^{-\beta J_{\beta}(\phi^2)+\beta\langle h,\phi^2\rangle},\qquad h\in\Ccal_{\rm b}(\R^d),\phi\in L^2(\R^d),\|\phi\|_2=1.
\end{equation}
If the random potential $h=q_N-W+f$ were in $\Ccal_{\rm b}(\R^d)$ almost surely, then \eqref{toestimate} instantly implied an estimate for the last term on the right hand side of \eqref{estimationstart}, and we could choose $\phi$ arbitrary and optimise later on $\phi$. 

However, the random potential $h=q_N-W+f$ does not have sufficient regularity for applying \eqref{toestimate} directly. But note that $q_j$ lies in $L^2(\R^d)$ almost surely, as is easily derived from the assumption that $\int_{\R^d}v(|x|)^2\,\d x<\infty$. In order to make \eqref{toestimate} applicable for functions $h$ of the form $q-W+f$ with $q\in L^2(\R^d)$, we have to establish a lower bound for $J_\beta$ in terms of a supremum over this class of potentials. This is the content of the following lemma.

\begin{lemma}\label{kato}For any $ \phi\in L^2(\R^d)\cap L^4(\R^d)$ satisfying $||\phi||_2=1 $, 
\begin{equation}\label{JRKato2}
J_\beta(\phi^2)\ge \sup_{h\in L^2(\R^d)\colon h\le 0}\Big(\langle -W+h,\phi^2\rangle-\frac{1}{\beta}\log\E\big[{\rm e}^{\beta\langle -W+h,\mu_\beta\rangle}\big]\Big).
\end{equation}
\end{lemma}

\begin{proofsect}{Proof} As a first step, we show that
\begin{equation}\label{JRKato}
J_{\beta}(\phi^2)\ge\sup_{f\in L^2(\R^d)\colon f\leq 0}\Big(\langle f,\phi^2 \rangle-\frac{1}{\beta}\log\E\big[{\rm e}^{\beta\langle f,\mu_\beta\rangle}\big]\Big).
\end{equation}
Let $ f\in  L^2(\R^d) $ be given satisfying $ f\le 0 $. We approximate $f$ by continuous bounded functions in a standard way as follows. Let $\kappa\colon\R^d\to[0,\infty) $ smooth with $ \int_{\R^d}\kappa(x)\,\d x=1 $ and $ \mbox{supp }\kappa\subset[-1,1]^d$. Fix $ \eps>0 $ and consider $ f_{\eps}=f\ast \kappa_{\eps} $, where $\kappa_{\eps}(x)=\eps^{-d}\kappa(x/\eps)$ for $x\in\R^d$. Then $f_\eps$ is continuous. Using Schwarz' inequality and the fact that $f\in L^2(\R^d)$, one sees that $f_\eps$ is bounded. Furthermore,
\begin{equation}\label{stineqKato}
\liminf_{\eps\to 0} \langle f_{\eps},\phi^2\rangle= \liminf_{\eps\to 0}\langle f,\phi^2\ast\kappa_{\eps}\rangle = \langle f,\phi^2\rangle,
\end{equation}
since $ f\in L^2(\R^d) $ and $ \phi^2\in L^2(\R^d) $. 

Since $f_\eps\to f$ strongly in $L^2$, we can pick a subsequence $\eps_n\downarrow 0$ such that $f_{\eps_n}\to f$ pointwise almost everywhere. Since $f\leq 0$, it follows from the bounded convergence theorem that
\begin{equation}\label{ndineqKato}
\lim_{n\to\infty}\E\big[{\rm e}^{\int_0^\beta f_{\eps_n}(B_s)\,\d s}\big]= \E\big[{\rm e}^{\int_0^\beta f(B_s)\,\d s}\big].
\end{equation}
From this, together with \eqref{stineqKato}, \eqref{JRKato} follows. 

Now we prove \eqref{JRKato2} by showing that
$$
\sup_{f\in L^2(\R^d)\colon f\leq 0}\Big(\langle f,\phi^2 \rangle-\frac{1}{\beta}\log\E\big[{\rm e}^{\beta\langle f,\mu_\beta\rangle}\big]\Big)
\geq \sup_{h\in L^2(\R^d)\colon h\le 0}\Big(\langle -W+h,\phi^2\rangle-\frac{1}{\beta}\log\E\big[{\rm e}^{\beta\langle -W+h,\mu_\beta\rangle}\big]\Big).
$$
This is similar to the proof of \lq$\geq$\rq\ in \eqref{JWbetacompare} in Lemma~\ref{Jbetaident}. Let $ h\in L^2(\R^d)$ satisfy $h\le 0 $, and consider $f_R:=(-W+h)\1_{Q_R} $ for $ R>0 $. Clearly, $f_R\in L^2(\R^d)$ with $f_R\leq 0$, and $f_R\downarrow -W+h$ pointwise as $R\to\infty$. Therefore, 
$$
\liminf_{R\to\infty}\langle f_R,\phi^2 \rangle\geq \langle -W+h,\phi^2\rangle,
$$
according to the monotonous convergence theorem. Furthermore, for the same reason, and since $-W+h$ is bounded from above,
$$
\limsup_{R\to\infty}\E\big[{\rm e}^{\beta\langle f_R,\mu_\beta\rangle}\big]
\leq \E\big[{\rm e}^{\beta\langle -W+h,\mu_\beta\rangle}\big].
$$
This implies the statement and finishes the proof of \eqref{JRKato2}.
\qed\end{proofsect}

We would like to remark that it is the assertion in \eqref{stineqKato} that forces us to require that $v\circ|\cdot|$ lies in $L^2(\R^d)$, since we will apply Lemma~\ref{kato} to the minimiser $\phi^2$ on the right-hand side of \eqref{chifdef}, and we do not know any higher integrability property of this function than that $\phi^2\in L^2(\R^d)$.

It is clear that \eqref{JRKato2} remains true if $W$ is replaced by $W-f$, where $f\in\Ccal_{\rm b}(\R^d)$, since adding a constant to $W-f$ does not change the value of the expression in the supremum on the right-hand side of \eqref{JRKato2}, and the potential $W-f-\inf f$ also satisfies Assumption (W).

Now we proceed with the proof of the lower bound of \eqref{logmomentfunctionKupper} in Proposition \ref{uppKN}. Go back to \eqref{estimationstart} and recall that the random potential $ q_j $ defined in (\ref{randompot}) is nonpositive and lies in $L^2(\R^d)$. Using Lemma~\ref{kato} with $ h=q_N $ and $ W $ replaced by $ W-f $, the last term on the right hand side of \eqref{estimationstart} is estimated as follows. For  any $ \phi\in L^2(\R^d)\cap L^4(\R^d) $ satisfying $ ||\phi||_2=1 $, 
\begin{equation}
\E^{\ssup{N}}\big[{\rm e}^{\beta\langle q_N-W+f,\mu_\beta^{\ssup{N}}\rangle}\big]\ge {\rm e}^{-\beta J_{\beta}(\phi^2)+\beta\langle q_N-W+f,\phi^2\rangle}.
\end{equation}  
Using this in (\ref{estimationstart}), we obtain
\begin{equation}\label{estimationstartb}
\begin{aligned}
\E\Big[&{\rm e}^{-H_{N,\beta}-K^{\ssup N}_{N,\beta}+\beta N\langle f,\overline{\mu}_{N,\beta}\rangle}\Big] 
\ge \E^{\ssup{1}}\otimes\cdots\otimes\E^{\ssup{N-2}}\Big[{\rm e}^{-H_{N-2,\beta}-K_{N-2,\beta}+\beta\langle  f,\mu_{N-2,\beta}\rangle}\\
&\E^{\ssup{N-1}}\big[{\rm e}^{\beta\langle q_{N-1}-W+f,\mu_\beta^{\ssup{N-1}}\rangle}\big]\Big]
\times {\rm e}^{-\beta J_{\beta}(\phi^2)-\beta\langle W-f,\phi^2\rangle-\beta\langle V_N\ast\phi^2,\phi^2\rangle}.
\end{aligned}
\end{equation}
Now we apply the same reasoning to the last expectation and iterate the argument. In this way we derive
\begin{equation}
\begin{aligned}
\frac{1}{N\beta}\log\E\Big[{\rm e}^{-H_{N,\beta}-K^{\ssup N}_{N,\beta}+\beta N\langle f,\overline{\mu}_{N,\beta}\rangle}\Big]
&\ge  -J_{\beta}(\phi^2)-\langle W-f,\phi^2\rangle -\frac{N-1}{2}\langle V_N\ast\phi^2,\phi^2\rangle.
\end{aligned}
\end{equation}
Since $ v\circ |\cdot |$ lies in $ L^1(\R^d) $ by assumption and since $ \phi^2\in L^2(\R^d) $, we have that 
$$
\lim_{N\to\infty}(N-1)V_N\ast\phi^2=\phi^2\int_{\R^d}v(|x|)\,\d x= \phi^2 4\pi\sca(v)\qquad \text{strongly in } L^2(\R^d).
$$
Hence,
\begin{equation}
\lim_{N\to\infty}\frac{N-1}{2}\langle V_N\ast\phi^2,\phi^2\rangle=4\pi\sca(v)\,||\phi||_4^4.
\end{equation} 
This implies that
\begin{equation}
\begin{aligned}
\liminf_{N\to\infty}\frac{1}{N\beta}\log\E\Big[{\rm e}^{-H_{N,\beta}-K^{\ssup N}_{N,\beta}+\beta N\langle f,\overline{\mu}_{N,\beta}\rangle}\Big] 
&\ge  - J_{\beta}(\phi^2)-\langle W-f,\phi^2\rangle-4\pi{\sca}(v)||\phi||_4^4\\ 
&\ge -\chi^{\ssup{\otimes}}_{{\sca}(v)}(f),
\end{aligned} 
\end{equation} 
and the proof of the lower bound in (\ref{logmomentfunctionKupper}) is finished.
\qed

\section{Appendix}
\subsection{Large deviations.}\label{LDP}

\noindent For the convenience of our reader, we repeat the notion of a large-deviation principle and of the most important facts that are used in the present paper. See \cite{DZ98} for a comprehensive treatment of this theory. 

Let $ \Xcal $ denote a topological vector space.  A lower semi-continuous function $ I\colon \Xcal\to [0,\infty] $ is called a {\it rate function\/} if  $ I $ is not identical $ \infty$ and  has compact level sets, i.e., if $ I^{-1}([0,c])=\{x\in\Xcal\colon I(x)\le c\} $ is compact for any $ c\ge 0 $. A sequence $(X_N)_{N\in\N}$ of $\Xcal$-valued random variables $X_N$  satisfies the {\it large-deviation upper bound\/} with {\it speed\/} $a_N$ and rate function $I$ if, for any closed subset $F$ of $\Xcal$,
\begin{equation}\label{LDPupper}
\limsup_{N\to\infty}\frac 1{a_N}\log \P(X_N\in F)\leq -\inf_{x\in F}I(x),
\end{equation}
and it satisfies the {\it large-deviation lower bound\/} if, for any open subset $G$ of $\Xcal$,
\begin{equation}\label{LDPlower}
\liminf_{N\to\infty}\frac 1{a_N}\log \P(X_N\in G)\leq -\inf_{x\in G}I(x).
\end{equation}
If both, upper and lower bound, are satisfied, one says that $(X_N)_N$ satisfies a {\it large-deviation principle}. The principle is called {\it weak\/} if the upper bound in \eqref{LDPupper} holds only for {\it compact\/} sets $F$. A weak principle can be strengthened to a full one by showing that the sequence of distributions of $X_N$ is {\it exponentially tight}, i.e., if for any $M>0$ there is a compact subset $K_M$ of $\Xcal$ such that $\P(X_N\in M^{\rm c})\leq {\rm e}^{-MN}$ for any $n\in\N$. 

One of the most important conclusions from a large deviation principle is {\it Varadhan's Lemma}, which says that, for any bounded and continuous function $F\colon \Xcal\to\R$,
$$
\lim_{N\to\infty}\frac 1N\log \int {\rm e}^{N F(X_N)}\,\d\P=-\inf_{x\in \Xcal}\big(I(x)-F(x)\big).
$$

All the above is usually stated for probability measures $\P$ only, but the notion easily extends to {\it sub}-probability measures $\P=\P_N$ depending on $N$. Indeed, first observe that the situation is not changed if $\P$ depends on $N$, since a large deviation principle depends only on distributions. Furthermore, the connection between probability distributions $\widetilde \P_N$ and sub-probability measures $\P_N$ is provided by the transformed measure $\widetilde \P_N(X_N\in A)=\P_N(X_N\in A)/\P_N(X_N\in\Xcal)$: If the measures $\P_N\circ X_N^{-1}$ satisfy a large deviation principle with rate function $I$, then the probability measures $\widetilde \P_N\circ X_N^{-1}$ satisfy a large deviation principle with rate function $I-\inf I$.

One standard situation in which a large deviation principle holds is the case where $\P$ is a probability measure, and $X_N=\frac 1N(Y_1+\dots+Y_N)$ is the mean of $N$ i.i.d.~$\Xcal$-valued random variables $Y_i$ whose moment generating function $M(F)=\int {\rm e}^{F(Y_1)}\,\d\P$ is finite for all elements $F$ of the topological dual space $\Xcal^*$ of $\Xcal$. In this case, the abstract {\it Cram\'er Theorem} provides a weak large deviation principle for $(X_N)_{N\in\N}$ with rate function equal to the Legendre-Fenchel transform of $\log M$, i.e., $I(x)=\sup_{F\in \Xcal^*}(F(x)-\log M(F))$. 

\subsection{Gross-Pitaevskii theory.}\label{GPtheory}

Consider the ground-state energy per particle of the Hamilton operator $ H_N$ in \eqref{HNdef},
\begin{equation}\label{chivdef}
\chi_{N }=\frac 1N\inf_{h\in H^1(\R^{dN})\colon\|h\|_2=1}\langle h,H_N h\rangle.
\end{equation}
It is standard to show the existence, uniqueness  and some regularity properties of the minimiser $h_N\in H^1(\R^{dN}) $. The large-$N$ behaviour, in a certain dilute regime, of $\chi_N$ and of the minimiser $h_N$ was studied by Lieb {\it et ~al.}~in a series of papers \cite{LSY00}, \cite{LY01}, \cite{LSY01}, \cite{LS02}, see also the monograph  \cite{LSSY05}. It turned out there that the Gross-Pitaevskii formula in \eqref{GP} well approximates the ground-state energy.   A summary of the large-$N$ results for $\chi_N$ is as follows. 
Assume that $d\in\{2,3\}$, that $v\geq 0$ with $v(0)\in(0,\infty]$, and $\int_{a+1}^\infty v(r)r^{d-1}\,\d r<\infty$, where $a=\inf\{r>0\colon v(r)<\infty\}\in[0,\infty)$. These assumptions guarantee that the {\em scattering length}, denoted by $ \widetilde{\sca}(v) $, is finite (\cite{LSSY05}). Note, that $ \sca(v) >\widetilde{\sca}(v) $ (\cite{ABK04}).

\begin{theorem}[Large-$N$ asymptotic of $\chi_{N }$ in $d\in\{2,3\}$, \cite{LSY00}, \cite{LY01}, \cite{LSY01}]\label{BECLieb}  Replace $v$ by $v_N(\cdot)=\beta_N^{-2}v(\,\cdot\,\beta_N^{-1})$ with $\beta_N=1/N$ in $d=3$ and $\beta_N^2=\widetilde \sca(v)^{-2}{\rm e}^{-N/\widetilde \sca(v)}N\|\phi^{\ssup{{\rm GP}}}_{\widetilde \sca(v)}\|_4^{-4}$ in $d=2$. Define $\phi^2_N\in H^1(\R^d)$ as the normalised  first marginal of $h_N^2$, i.e.,
$$
\phi_N^2(x)=\int_{\R^{d(N-1)}}h_N^2(x,x_2,\dots ,x_N)\,\d x_2\cdots\d x_N,\qquad x\in\R^d.
$$
Then
$$
\lim_{N\to\infty}\chi_{N }=\chi^{\ssup{\rm GP}}_{\widetilde \sca(v)}\qquad\mbox{and}\qquad \phi_N^2\to\big(\phi^{\ssup{{\rm GP}}}_{\widetilde \sca(v)})^2\quad\mbox{ in weak $L^1(\R^d)$-sense.}
$$
\end{theorem}

The proofs show that the ground state, $h_N$, approaches the product-state $ (\phi^{\ssup{\rm GP}}_{\widetilde \sca(v)})^{\otimes N}$ if $N$ gets large.

Moreover, on the basis of this result, the occurrence of Bose-Einstein condensation in the ground-state (zero-temperature) was proved in \cite{LS02}.

\subsection*{Acknowledgments} 

\noindent This work was partially supported by DFG grant AD 194/1-1.

\end{document}